\documentclass[english]{article}
\usepackage[T1]{fontenc}
\usepackage[latin9]{inputenc}
\usepackage{geometry}
\usepackage[active]{srcltx}
\usepackage{float}
\usepackage{textcomp}
\usepackage{amsmath}
\usepackage{amsthm}
\usepackage{amssymb}
\usepackage{graphicx}
\usepackage{color}

\makeatletter

\providecommand{\tabularnewline}{\\}
\floatstyle{ruled}
\newfloat{algorithm}{tbp}{loa}
\providecommand{\algorithmname}{Algorithm}
\floatname{algorithm}{\protect\algorithmname}

\numberwithin{equation}{section}
\numberwithin{figure}{section}

\usepackage{algorithm,algpseudocode}
\allowdisplaybreaks

\makeatother

\usepackage{babel}
\begin{document}
\global\long\def\ve{\varepsilon}
\global\long\def\R{\mathbb{R}}
\global\long\def\Rn{\mathbb{R}^{n}}
\global\long\def\Rd{\mathbb{R}^{d}}
\global\long\def\E{\mathbb{E}}
\global\long\def\P{\mathbb{P}}
\global\long\def\bx{\mathbf{x}}
\global\long\def\vp{\varphi}
\global\long\def\ra{\rightarrow}
\global\long\def\smooth{C^{\infty}}
\global\long\def\Tr{\mathrm{Tr}}
\global\long\def\bra#1{\left\langle #1\right|}
\global\long\def\ket#1{\left|#1\right\rangle }

\newcommand{\ML}[1]{\textcolor{red}{[ML:#1]}}
\newcommand{\YK}[1]{\textcolor{blue}{[YK:#1]}}

\title{Scalable semidefinite programming approach to \\
variational embedding for quantum many-body problems }

\author{Yuehaw Khoo and Michael Lindsey}
\maketitle
\begin{abstract}
In quantum embedding theories, a quantum many-body system
is divided into localized clusters of sites which are treated with
an accurate `high-level' theory and glued together self-consistently
by a less accurate `low-level' theory at the global scale. The recently
introduced variational embedding approach for quantum many-body problems
combines the insights of semidefinite relaxation and quantum embedding
theory to provide a lower bound on the ground-state energy that improves
as the cluster size is increased. The variational embedding method
is formulated as a semidefinite program (SDP), which can suffer from
poor computational scaling when treated with black-box solvers. We
exploit the interpretation of this SDP as an embedding method to develop
an algorithm which alternates parallelizable local updates of the
high-level quantities with updates that enforce the low-level global
constraints. Moreover, we show how translation invariance in lattice systems can be exploited to reduce the complexity of projecting a key matrix to the positive semidefinite cone.
\end{abstract}

\section{Introduction}

The problem of determining the ground state of a quantum many-body
system has wide-ranging applications in physics, chemistry, and materials
science. This problem can be viewed as the problem of determinining
the lowest eigenvalue of a Hermitian operator on a Hilbert space whose
dimension grows exponentially with the size of the system or the number
of particles, such as electrons in the case of electronic structure.
Here we highlight two relevant categories of approaches to taming
this curse of dimensionality.

The first category is that of semidefinite relaxations, which rephrase
the aforementioned energy minimization problem as an optimization
problem in terms of a reduced set of physical observables, almost
always a semidefinite program (SDP). In principle these observables
satisfy representability constraints, i.e., constraints that ensure
that they can be recovered from a bona fide quantum many-body state.
However, only a subset of representability constraints can be efficiently
enforced, yielding tractable optimization problems that provide lower
bounds on the ground-state energy. Such approaches include the 2-RDM
theories~\cite{Mazziotti1998,Mazziotti2004,CancesStoltzLewin2006,Mazziotti2012,ZhaoBraamsFukudaEtAl2004,LiWenYangEtAl2018,AndersonNakataIgarashiEtAl2013,NakataNakatsujiEharaEtAl2001,DeprinceMazziotti2010}, as well as methods that may be classified as
quantum marginal relaxations such as~\cite{LeiferPoulin2008,PoulinHastings2011,FerrisPoulin2013,VariationalCorrelations}.

Meanwhile, quantum embedding theories take the perspective of dividing
a system into local clusters, small enough to be treated with a highly
accurate or exact method referred to as the `high-level' method. Local
problems are then stitched together via a reduced set of global quantities
or a less-accurate `low-level' method that operates on the global
scale, and the local and global perspectives are constrained to be
compatible via some self-consistency condition. Such approaches include
dynamical mean-field theory (DMFT)~\cite{GeorgesKotliarKrauthEtAl1996,KotliarSavrasovHauleEtAl2006} and density matrix embedding
theory (DMET)~\cite{KniziaChan2012,KniziaChan2013}, as well as variants such as the energy-weighted
DMET (EwDMET)~\cite{FertittaBooth2018,FertittaBooth2019} which in a certain sense interpolates between
DMFT and DMET~\cite{SriluckshmyEtAl2020}.

Recently, variational embedding~\cite{LinLindsey2020} was introduced as a semidefinite
relaxation which is also a quantum embedding method. Like other relaxations
such as~\cite{LeiferPoulin2008,PoulinHastings2011,FerrisPoulin2013,VariationalCorrelations}, the key optimization variables are quantum marginals
for local clusters, but variational embedding additionally includes
\emph{global }constraints which tighten the relaxation and accommodate
the treatment of, e.g., long-range interactions.

\subsection{Contribution}

As an SDP, variational embedding can be solved with black-box methods,
as is done in~\cite{LinLindsey2020}, but scalability calls for a solver that
is specially adapted to the problem. In this work, we introduce a
scalable solver for this SDP which takes advantage of the embedding
interpretation of the approach. The aforementioned global constraint
is dualized to reduce the problem to a simpler relaxation (similar
to those of~\cite{LeiferPoulin2008,PoulinHastings2011,FerrisPoulin2013,VariationalCorrelations}) in which the global constraints have been
exchanged for effective contributions to relevant Hamiltonian operators
at the local scale. This problem can then be solved in a fashion in
which effective problems for the key variables (the two-cluster marginals)
are completely decoupled and can be solved in parallel, with dual
variable optimization enforcing the self-consistency of these problems. Furthermore, translation invariance of a lattice system can be used to significantly speed up the running time.

Our approach is based on augmented Lagrangian methods (see for example \cite{boyd2011distributed,sun2015convergent,sun2020sdpnal}), which have been used for solving large-scale SDP problems. In particular, ADMM-type approaches can allow for sub-problems to be solved in parallel. Approaches such as \cite{sun2015convergent} apply ADMM to the dual problem, whereas in \cite{wen2010alternating} the primal problem is solved. Our approach differs in that the local constraints are kept in their primal form while the global positive semidefinite constraint that couples the local variables is taken into account via the introduction a dual variable. Consequently, each iteration involves the solution of many decoupled effective problems, preserving the flavor of a quantum embedding theory.

\subsection{Outline}

In Section \ref{sec:Preliminaries} we provide relevant background
on the ground-state eigenvalue problem, examples of interest, and
the two-marginal relaxation for variational embedding introduced in~\cite{LinLindsey2020}. (In Appendix \ref{sec:Fermions}, additional background
is provided for the context of fermionic systems.) In Section \ref{sec:Optimization-approach},
we describe our optimization approach to this problem, which is an
SDP. The section begins with an idealized scheme of projected gradient
ascent on the dual variable to the aforementioned global constraint.
In order to implement such a scheme, it is necessary to solve an effective
problem in terms of the primal variables. In Section \ref{sub:ADMM-type-approach-to},
we introduce an ADMM-type approach to this problem, and in Section
\ref{sub:Practical-scheme} we integrate this approach with dual
ascent to define our practical scheme. In Section \ref{sub:Exploiting-translation-invarianc}
we explain how translation-invariance can be exploited, and in Section
\ref{sub:Discussion-of-scaling} we include a detailed discussion
of the computational scaling. Finally in Section \ref{sec:Numerical-experiments}
we present numerical experiments on several model systems of quantum
spins and fermions.

\subsection*{Acknowledgments}
This work was partially supported by the National Science Foundation under Award No. 1903031 (M.L.). We thank Lin Lin for helpful discussions.

\section{Preliminaries\label{sec:Preliminaries}}

In this section we review the formulation of variational embedding
for quantum spins, following~\cite{LinLindsey2020}. In Appendix \ref{sec:Fermions},
we review the case of fermions (also following~\cite{LinLindsey2020}), which
requires a bit more care but nonetheless yields a semidefinite program
of identical form after suitable manipulations.

\subsection{The ground-state eigenvalue problem}

We consider a model with $M$ sites, indexed $i=1,\ldots,M$, each
endowed with a classical local state space $X_{i}$ (which shall be
discrete). These in turn yields local quantum state spaces $Q_{i}=\mathbb{C}^{\vert X_{i} \vert }$.
The global quantum state space (i.e., the space of wavefunctions)
is then given by 
\[
\mathcal{Q}:=\bigotimes_{i=1}^{M}Q_{i}\simeq\mathbb{C}^{\vert \mathcal{X} \vert },
\]
 where $\mathcal{X}:=\prod_{i=1}^{M}X_{i}$ is the global classical
state space, so quantum states (wavefunctions) correspond to complex-valued
functions on the classical state space. Let $H_{i}$ (resp. $H_{ij}$)
denote Hermitian operators $Q_{i}\ra Q_{i}$ (resp. $Q_{ij}\ra Q_{ij}$),
and let $\hat{H}_{i}$ (resp. $\hat{H}_{ij}$) denote the corresponding
operators $\mathcal{Q}\ra\mathcal{Q}$ obtained by tensoring $H_{i}$
by the identity operator on all sites $k\neq i$ (resp. $k\neq i,j$).
We consider \emph{pairwise} Hamiltonian $\hat{H}:\mathcal{Q}\ra\mathcal{Q}$
of the form 
\[
\hat{H}=\sum_{i}\hat{H}_{i}+\sum_{i<j}\hat{H}_{ij},
\]
 and our interest is in determining the ground-state energy, i.e.,
the lowest eigenvalue, of $\hat{H}$. We denote this eigenvalue by
$E_{0}$, which is defined variationally by 
\begin{equation}
E_{0}=\inf\left\{ \Phi^{*}\hat{H}\Phi\,:\,\Phi\in\mathcal{Q},\:\Phi^{*}\Phi=1\right\} .\label{eq:groundstate}
\end{equation}

Now we review some examples of interest. First consider the case of
quantum spin-$\frac{1}{2}$ systems, i.e., the case $X_{i}=\{-1,1\}$.
To construct operators on $\mathcal{Q}$, one first starts with the
Pauli matrices 
\[
\sigma^{x}=\left(\begin{array}{cc}
0 & 1\\
1 & 0
\end{array}\right),\quad\sigma^{y}=\left(\begin{array}{cc}
0 & -i\\
i & 0
\end{array}\right),\quad\sigma^{z}=\left(\begin{array}{cc}
1 & 0\\
0 & -1
\end{array}\right),
\]
 which (together with the identity $I_{2}$) form a basis for the
real vector space of Hermitian operators on $\mathbb{C}^{2}$. We
let $\sigma_{i}^{x/y/z}$ denote the operator $\mathcal{Q}\ra\mathcal{Q}$
obtained by tensoring $\sigma^{x/y/z}$ on the $i$-th site with $I_{2}$
on all other sites. Then in terms of these operators we can define
the transverse-field Ising (TFI) Hamiltonian and the anti-ferromagnetic
Heisenberg (AFH) Hamiltonian by 
\begin{equation}
\hat{H}_{\mathrm{TFI}}=-h\sum_{i}\sigma_{i}^{x}-\sum_{i\sim j}\sigma_{i}^{z}\sigma_{j}^{z}\label{eq:TFI}
\end{equation}
 
\begin{equation}
\hat{H}_{\mathrm{AFH}}=\sum_{i\sim j}\left[\sigma_{i}^{x}\sigma_{j}^{x}+\sigma_{i}^{y}\sigma_{j}^{y}+\sigma_{i}^{z}\sigma_{j}^{z}\right],\label{eq:AFH}
\end{equation}
 $h$ is a scalar paramater and summation over $i\sim j$ indicates
summation over pairs of indices that are adjacent within some graph
defined on the index set $\{1,\ldots,M\}$, often a rectangular lattice
in some dimension. These problems have been considered as prototypical quantum many-body problems, e.g., in \cite{CarleoTroyer2017}, as well as models for the study of quantum phase transitions, as in \cite{Sachdev}.

\subsection{The two-marginal relaxation}

In~\cite{LinLindsey2020}, the optimization problem (\ref{eq:groundstate}) is reformulated
as an optimization over the density operator $\rho$, which (for nondegenerate
ground states) corresponds at optmality to $\Phi_{0}\Phi_{0}^{*}$
where $\Phi_{0}$ is the ground state eigenvector, i.e., the optimizer
of (\ref{eq:groundstate}). This optimization is in turn relaxed as
a computationally tractable optimization over the quantum two-marginals
\[
\rho_{ij}=\Tr_{\{i,j\}^{c}}[\rho],\quad i<j,
\]
 which are defined as partial traces of $\rho$, analogous to classical
marginals. Indeed, recall~\cite{LinLindsey2020} that for any subset $S\subset\{1,\ldots,M\}$,
the partial trace $\rho_{S}=\Tr_{S^{c}}[\rho]$ may be defined as
the unique operator on $\bigotimes_{i\in S}Q_{i}$ such that $\Tr[A\rho_{S}]=\Tr[\hat{A}\rho]$
for all operators $A$ on $\bigotimes_{i\in S}Q_{i}$ (lifted to operators
$\hat{A}$ on $\mathcal{Q}$ by tensoring with the identity on $S^{c}$).
In particular, $\rho_{ij}$ is an operator on $Q_{i}\otimes Q_{j}$,
and moreover it is positive semidefinite with unit trace (following
from the same properties for $\rho$).\\

Then the two-marginal relaxation of~\cite{LinLindsey2020} reads in terms of the
two marginals (and the analogously-defined one-marginals, which can
be obtained from the two-marginals by further partial trace) as the
following semidefinite program, whose optimal value we denote by $E_{0}^{(2)}$:
\begin{eqnarray}
\underset{\{\rho_{i}\},\,\{\rho_{ij}\}_{i<j}}{\mbox{minimize}} &  & \sum_{i}\Tr\left[H_{i}\rho_{i}\right]+\sum_{i<j}\Tr\left[H_{ij}\rho_{ij}\right]\label{eq:sdpObj}\\
\mbox{subject to} &  & \rho_{ij}\succeq0,\quad1\leq i<j\leq M,\label{eq:sdpC1}\\
 &  & \rho_{i}=\Tr_{\{2\}}[\rho_{ij}],\ \rho_{j}=\Tr_{\{1\}}[\rho_{ij}],\quad1\leq i<j\leq M,\label{eq:sdpC2}\\
 &  & \Tr[\rho_{i}]=1,\quad i=1,\ldots,M,\label{eq:sdpC3}\\
 &  & G[\{\rho_{i}\},\{\rho_{ij}\}_{i\leq j}]\succeq0.\label{eq:sdpC4}
\end{eqnarray}
 Here $G=G[\{\rho_{i}\},\{\rho_{ij}\}_{i\leq j}]$ is an operator
defined linearly in terms of the one- and two-marginals, subordinate
to the specification of an arbitrary collection $\{O_{i,\alpha}\,:\,\alpha=1,\ldots,n_{i}\}$
of linear operators $Q_{i}\ra Q_{i}$ at each site $i=1,\ldots,M$.
In specific, $G$ is specified blockwise, with blocks $G_{ij}$ for
$1\leq i,j,\leq M$ of size $n_{i}\times n_{j}$ defined by 
\[
\left(G_{ij}\right)_{\alpha\beta}=\begin{cases}
\Tr\left[O_{i,\alpha}^{\dagger}O_{i,\beta}\,\rho_{i}\right] & i=j\\
\Tr\left[\left(O_{i,\alpha}^{\dagger}\otimes O_{j,\beta}\right)\rho_{ij}\right] & i\neq j.
\end{cases}
\]
 The choice of operators only matters up to $\mathrm{span}\{O_{i,\alpha}\,:\,\alpha=1,\ldots,n_{i}\}$,
and in our numerical experiments we shall consider the complete operator
collection spanning all linear maps $Q_{i} \ra Q_{i}$. The last constraint
(\ref{eq:sdpC4}) is called the \emph{global semidefinite constraint}.

\subsubsection{Classical marginal relaxation}
To motivate the relaxation \eqref{eq:sdpObj} further, we examine the problem of finding the lowest energy state of a classical energy function of a pairwise form. (This can be viewed as a special case of the more general quantum ground-state problem by taking the $H_{ij}$ to be diagonal operators.) More concretely, for $x_i \in X_i,\,i=1,\ldots,M$, define an energy function  
\begin{equation}
    E(x_1,\ldots,x_M) = \sum_{i<j} E_{ij}(x_i,x_j),
\end{equation}
and observe that the minimizer of $E$ can be determined via the linear program 
\begin{equation}
    \underset{\mu\in \Pi(\mathcal{X})}{\arg\min} \  \sum_{(x_1,\ldots,x_M) \in \mathcal{X}} E(x_1,\ldots,x_M) \, \mu(x_1,\ldots,x_M), 
\end{equation}
where $\Pi (\mathcal{X})$ is the space of probability measures on $\mathcal{X}$. Indeed, the optimizer is a $\delta$-function supported on the minimizer of $E$ (provided that it is unique). Exploiting the pairwise structure of $E$, we have 
\begin{equation}
    \underset{\{\mu_{ij}\}_{i<j} \,\mathrm{repr.}}{\arg\min} \ \sum_{i<j} \sum_{x_i,x_j\in \mathcal{X}} E_{ij}(x_i,x_j) \, \mu_{ij}(x_i,x_j)
\end{equation}
where the two-marginal variables $\{\mu_{ij}\}_{i<j}$ are constrained to be jointly representable, i.e., to be derivable as the two-marginals of a high-dimensional measure $\mu$. Enforcing this constraint demands exponential complexity, so various convex relaxation approaches have been proposed, where only certain necessary conditions for the $\{\mu_{ij}\}_{i<j}$ are kept; see, for instance, \cite{WainwrightJordan} for a review. In particular, the analogous convex relaxation to \eqref{eq:sdpObj} is 
\begin{eqnarray}
\underset{\{\mu_{i}\},\,\{\mu_{ij}\}_{i<j}}{\mbox{minimize}} &  & \sum_{i<j}\Tr\left[E_{ij}\mu_{ij}\right]\label{eq:csdpObj}\\
\mbox{subject to} &  & \mu_{ij}\geq 0,\quad1\leq i<j\leq M,\label{eq:csdpC1}\\
 &  & \mu_{i}=\mu_{ij}\mathbf{1},\ \mu_{j}=\mu_{ij}^\top \mathbf{1},\quad1\leq i<j\leq M,\label{eq:csdpC2}\\
 &  & \mu_{i}^\top \mathbf{1}=1,\quad i=1,\ldots,M,\label{eq:csdpC3}\\
 &  & \begin{bmatrix}\text{diag}(\mu_1) & \mu_{12} &\cdots & \mu_{1M} \\ \mu_{21} &  & & \\ \vdots & & \ddots& \vdots \\ \mu_{M1} & \cdots & & \text{diag}(\mu_M) \end{bmatrix}\succeq0, \label{eq:csdpC4}
\end{eqnarray}
where $\mathbf{1}$ is an all-one vector of appropriate size. Here \eqref{eq:csdpC1} and \eqref{eq:csdpC3} are standard constraints for any discrete probability distributions, and \eqref{eq:csdpC2} are constraints that enforces `local consistency' \cite{WainwrightJordan} among the $\{\mu_{ij}\}_{i<j}$. The global semidefinite constraint \eqref{eq:csdpC4} is discussed in \cite{KhooEtAl2019} and \cite{ChenEtAl2020} in the contexts of multi-marginal optimal transport and energy minimization, respectively.

\subsection{Cluster relaxation}

Given a quantum spin model as above and a decomposition of the sites
$\{1,\ldots,M\}$ as a disjoint union $\bigcup_{i'=1}^{M'}C_{i'}$
of clusters $C_{i'}$, we may define $X_{i'}':=\prod_{i\in C_{i'}}C_{i'}$
to be the classical state space for the $i'$-th cluster. One see
that any Hamiltonian that is pairwise with respect to sites is pairwise
with respect to clusters, so by viewing our clusters as sites and
applying the above formalism, we obtain a tighter relaxation~\cite{LinLindsey2020}.

\subsection{Partial duality}

In~\cite{LinLindsey2020} it was shown that (\ref{eq:sdpObj}) admits the minimax
formalization (obtained via dualization of the global semidefinite
constraint (\ref{eq:sdpC4})) 
\begin{equation}
E_{0}^{(2)}=\sup_{X\succeq0}\ \mathcal{F}[X],\label{eq:partialDual}
\end{equation}
 where 
\begin{equation}
\mathcal{F}[X]:=\inf\left\{ \sum_{i}\Tr\left(H_{i}[X_{ii}]\rho_{i}\right)+\sum_{i<j}\Tr\left(H_{ij}[X_{ij}]\rho_{ij}\right)\,:\,\{\rho_{i}\},\,\{\rho_{ij}\}_{i<j}\ \mbox{satisfy}\ \mbox{\eqref{eq:sdpC1}-\eqref{eq:sdpC3}}\right\} .\label{eq:partialDualObj}
\end{equation}
 Here $X_{ij}$ denote the blocks of $X$, and the `effective' Hamiltonian
terms $H_{i}[X_{ii}]$ and $H_{ij}[X_{ij}]$ are defined linearly
in terms of $X$ via 
\begin{equation}
H_{i}[X_{ii}]:=H_{i}-\sum_{\alpha,\beta=1}^{n_{i}}(\overline{X}_{ii})_{\alpha\beta}O_{i,\alpha}^{\dagger}O_{i,\beta},\quad H_{ij}[X_{ij}]:=H_{i}-\left[\sum_{\alpha=1}^{n_{i}}\sum_{\beta=1}^{n_{j}}(\overline{X}_{ij})_{\alpha\beta}\left(O_{i,\alpha}^{\dagger}\otimes O_{i,\beta}\right)+\mathrm{h.c.}\right],\label{eq:effHam}
\end{equation}
 where `h.c.' denotes the Hermitian conjugate term.
 
This partial dual formulation can be obtained from \eqref{eq:sdpObj} by exchanging the global semidefinite constraint \eqref{eq:sdpC4} for an extra term
\[
- \Tr[X\, G[\{\rho_{i}\},\{\rho_{ij}\}_{i\leq j}]]
\]
in the Lagrangian, where $X \succeq 0 $ is a dual variable with respect to which the Lagrangian is to be maximized. Then \eqref{eq:effHam} is recovered by breaking this additional term into a blockwise sum, collecting terms, and minimizing over the primal variables, subject to the remaining constraints.

 Since $X$ and
$G$ are dual variables~\cite{LinLindsey2020} we have that 
\begin{equation}
\nabla_{X}\mathcal{F}[X]=-G[\{\rho_{i}\},\{\rho_{ij}\}_{i<j}],\label{eq:gradX}
\end{equation}
 where $\{\rho_{i}\},\{\rho_{ij}\}_{i<j}$ are the minimizers of the
infimum in (\ref{eq:partialDualObj}).

\section{Optimization approach\label{sec:Optimization-approach}}

In order to solve the two-marginal relaxation (\ref{eq:sdpObj}),
our point of departure is the partial dual formulation (\ref{eq:partialDual}).
For simplicity we consider the case in which the Hamiltonian and all
variables are purely real, and we let $\Pi_{\succeq0}$
be defined by
\begin{equation}
    \Pi_{\succeq 0}(C) := \min_{S\succeq 0} \|S - C \|_F^2,
\end{equation}
for a symmetric matrix $A$, i.e. the Euclidean projection (in Frobenius norm) of a symmetric matrix onto
the set of real symmetric positive semidefinite matrices (equivalent
to setting all negative eigenvalues of the argument to zero). As an
idealized scheme, we can imagine performing projected gradient ascent
on (\ref{eq:partialDual}), which is implemented by Algorithm \ref{alg:gradAscent}. 

\begin{algorithm}
\caption{Exact projected gradient ascent}

\begin{algorithmic}[1]
\Require{$\ve > 0,\  X \succeq 0$}
\While{not converged}
\State{Set $\left( \{\rho_{i}\},\{\rho_{ij}\}_{i<j} \right)$ to be the minimizer in (\ref{eq:partialDualObj}), holding $X$ fixed}\label{eq:exactBP}
\State{$X\leftarrow\Pi_{\succeq0}\left(X+\ve G[\{\rho_{i}\},\{\rho_{ij}\}_{i<j}]\right)$}
\EndWhile
\end{algorithmic}\label{alg:gradAscent}
\end{algorithm}

The focus of this section is in the development of algorithm for step \ref{eq:exactBP} for a general Hamiltonian, and we also study the case in the presence of translational invariance. In practice, we will not fully converge a solution to step \ref{eq:exactBP}
of Algorithm \ref{alg:gradAscent}, resulting in an inexact projected
gradient ascent scheme. However, in order to motivate our practical
scheme, we will first discuss how to solve step \ref{eq:exactBP}
exactly for fixed $X$.

\subsection{Details for step 2 in algorithm 1\label{sub:ADMM-type-approach-to}}

We rephrase step \ref{eq:exactBP} as the following optimization problem:
\begin{eqnarray}
\underset{\{\rho_{i}\},\,\{\rho_{ij}\}_{i<j}}{\mbox{minimize}} &  & \sum_{i}\Tr\left[H_{i}'\rho_{i}\right]+\sum_{i<j}\Tr\left[H_{ij}'\rho_{ij}\right]\label{eq:BP}\\
\mbox{subject to} &  & \rho_{ij}\succeq0,\quad1\leq i<j\leq M,\nonumber \\
 &  & \rho_{i}=A_{1}[\rho_{ij}],\ \rho_{j}=A_{2}[\rho_{ij}],\quad1\leq i<j\leq M,\nonumber \\
 &  & \Tr[\rho_{i}]=1,\quad i=1,\ldots,M,\nonumber 
\end{eqnarray}
 where $H_{i}':=H_{i}[X_{ii}]$, $H_{ij}':=H_{ij}[X_{ij}]$, and $X$
is fixed for the duration of this subsection. Moreover, for simplicity
we have assumed that $m:=\vert X_{i}\vert$ is constant, and $A_{1},A_{2}$
are defined to be the linear operators $\Tr_{\{2\}}$ and $\Tr_{\{1\}}$,
respectively. Hence $A_{1},A_{2}$ can be realized as sparse matrices
of size $m^{2}\times m^{4}$ We then formulate an equivalent optimization
problem via the introduction of dummy variables $\tilde{\rho}_{ij}$
and the inclusion of augmented Lagrangian terms in the objective:
\begin{eqnarray}
\underset{\{\rho_{i}\},\,\{\rho_{ij},\tilde{\rho}_{ij}\}_{i<j}}{\mbox{minimize}} &  & \sum_{i}\Tr\left[H_{i}'\rho_{i}\right]+\sum_{i<j}\Tr\left[H_{ij}'\rho_{ij}\right]\nonumber \\
 &  & \quad+\ \sum_{i<j}\left(\frac{\mu}{2}\Vert\rho_{ij}-\tilde{\rho}_{ij}\Vert_{\mathrm{F}}^{2}+\frac{\nu}{2}\Vert\rho_{i}-A_{1}[\rho_{ij}]\Vert_{\mathrm{F}}^{2}+\frac{\nu}{2}\Vert\rho_{j}-A_{2}[\rho_{ij}]\Vert_{\mathrm{F}}^{2}\right)\label{eq:sdpALMobj}\\
\mbox{subject to} &  &   \tilde{\rho}_{ij}\succeq0,\quad1\leq i<j\leq M,\nonumber \\
 &  &  \Lambda_{ij}\ :\ \rho_{ij}=\tilde{\rho}_{ij},\quad1\leq i<j\leq M,\label{eq:dummy}\\
 &  & \Lambda_{ij}^{(1)}:\ \rho_{i}=A_{1}[\rho_{ij}],\ \Lambda_{ij}^{(2)}:\ \rho_{j}=A_{2}[\rho_{ij}],\quad1\leq i<j\leq M,\label{eq:lcc}\\
 &  & \Tr[\rho_{i}]=1,\quad i=1,\ldots,M,\nonumber 
\end{eqnarray}
 where $\mu,\nu>0$ are constant parameters, and  $\Lambda_{ij},\Lambda_{ij}^{(1)},\Lambda_{ij}^{(2)}\in\mathrm{End}(Q_{i}\otimes Q_{j})$
are the dual variables for the associated constraints in (\ref{eq:dummy}) and (\ref{eq:lcc}). Let $f[\{\rho_{i}\},\,\{\rho_{ij},\tilde{\rho}_{ij}\}_{i<j}]$
denote the objective function (\ref{eq:sdpALMobj}), yielding the Lagrangian 
\begin{equation}
\begin{aligned} & \mathcal{L}\left(\{\rho_{i}\},\,\{\rho_{ij},\tilde{\rho}_{ij}\}_{i<j};\,\{\Lambda_{ij},\Lambda_{ij}^{(1)},\Lambda_{ij}^{(2)}\}_{i<j}\right)\\
 & \ \ :=\ f[\{\rho_{i}\},\,\{\rho_{ij},\tilde{\rho}_{ij}\}_{i<j}]+\sum_{i<j}\left(\left\langle \Lambda_{ij},\tilde{\rho}_{ij}-\rho_{ij}\right\rangle _{\mathrm{F}}+\left\langle \Lambda_{ij}^{(1)},A_{1}[\rho_{ij}]-\rho_{i}\right\rangle _{\mathrm{F}}+\left\langle \Lambda_{ij}^{(2)},A_{2}[\rho_{ij}]-\rho_{j}\right\rangle _{\mathrm{F}}\right),
\end{aligned}
\label{eq:lagrangian}
\end{equation}
with domain specified by the (undualized) primal constraints $\tilde{\rho}_{ij}\succeq0$
and $\Tr[\rho_{i}]=1$. Here $\left\langle \,\cdot\,,\,\cdot\,\right\rangle _{\mathrm{F}}$
indicates the Frobenius inner product. Then the Augmented Lagrangian method \cite{bertsekas2014constrained} for (\ref{eq:BP}) is implemented by Algorithm \ref{alg:ALM}.

\begin{algorithm}
\caption{Augmented Lagrangian method for (\ref{eq:BP})}

\begin{algorithmic}[1]
\Require{$\mu,\nu > 0, \ \{H'_i\},\  \{ H'_{ij}, \Lambda_{ij}, \Lambda_{ij}^{(1)}, \Lambda_{ij}^{(2)}\}_{i<j}$}
\While{not converged}
%\State{Set $\left(\{\rho_{i}\},\{\rho_{ij},\tilde{\rho}_{ij}\}_{i<j}\right)$ to be the minimizer of \eqref{eq:lagrangian}, subject to the constraints that $\Tr[\rho_{i}]=1$ and $\tilde{\rho}_{ij}\succeq0$ for all $i<j$, holding  $\{\Lambda_{ij},\Lambda_{ij}^{(1)},\Lambda_{ij}^{(2)}\}_{i<j}$ fixed} \label{eq:exactALM}
\State{$\{\rho_{i}\},\{\rho_{ij},\tilde{\rho}_{ij}\}_{i<j}\leftarrow \underset{{\{\rho_{i}\},\{\rho_{ij},\tilde{\rho}_{ij}\}_{i<j}}}{\arg\min} \mathcal{L}\left(\{\rho_{i}\},\,\{\rho_{ij},\tilde{\rho}_{ij}\}_{i<j};\,\{\Lambda_{ij},\Lambda_{ij}^{(1)},\Lambda_{ij}^{(2)}\}_{i<j}\right)$ } \label{eq:exactALM}
\For{each pair $i<j$}
\State{$\Lambda_{ij}\leftarrow\Lambda_{ij}+\mu\left(\tilde{\rho}_{ij}-\rho_{ij}\right)$
}
\State{$\Lambda_{ij}^{(1)}\leftarrow\Lambda_{ij}^{(1)}+\nu(A_{1}[\rho_{ij}]-\rho_{i})$
}
\State{$\Lambda_{ij}^{(2)}\leftarrow\Lambda_{ij}^{(2)}+\nu(A_{2}[\rho_{ij}]-\rho_{j})$
}
\EndFor

\EndWhile
\end{algorithmic}

\label{alg:ALM}
\end{algorithm}

In practice, it is difficult to solve step \ref{eq:exactALM} of Algorithm
\ref{alg:ALM} exactly. Therefore, instead of optimizing $\{\rho_{i}\},\{\rho_{ij},\tilde{\rho}_{ij}\}_{i<j}$
jointly, we consider an ADMM-type \cite{boyd2011distributed} substitute, namely Algorithm \ref{alg:ADMM}.
\begin{algorithm}
\caption{Pseudo-code for ADMM-type method for (\ref{eq:BP})}

\begin{algorithmic}[1]
\Require{$\mu,\nu > 0, \ \{H'_i, \rho_i\},\ \{H'_{ij}, \tilde{\rho}_{ij}, \Lambda_{ij}, \Lambda_{ij}^{(1)}, \Lambda_{ij}^{(2)}\}_{i<j}$}
\While{not converged}
\State{$\{\rho_{ij}\}_{i<j}\leftarrow \underset{{\{\rho_{ij}\}_{i<j}}}{\arg\min}\  \mathcal{L}\left(\{\rho_{i}\},\,\{\rho_{ij},\tilde{\rho}_{ij}\}_{i<j};\,\{\Lambda_{ij},\Lambda_{ij}^{(1)},\Lambda_{ij}^{(2)}\}_{i<j}\right)$ } \label{eq:ADMM1}
\State{$\{\rho_{i}\},\{\tilde{\rho}_{ij}\}_{i<j}\leftarrow \underset{{\{\rho_{i}\},\{\tilde{\rho}_{ij}\}_{i<j}}}{\arg\min} \mathcal{L}\left(\{\rho_{i}\},\,\{\rho_{ij},\tilde{\rho}_{ij}\}_{i<j};\,\{\Lambda_{ij},\Lambda_{ij}^{(1)},\Lambda_{ij}^{(2)}\}_{i<j}\right)$}\label{eq:ADMM2}
%\State{Set $\{\rho_{ij} \}_{i<j}$ to be the minimizer of \eqref{eq:lagrangian}, holding  $\{\rho_{i}\}$ and $\{\tilde{\rho}_{ij},\Lambda_{ij},\Lambda_{ij}^{(1)},\Lambda_{ij}^{(2)}\}_{i<j}$ fixed} \label{eq:ADMM1}
%\State{Set $\left(\{\rho_{i}\},\{\tilde{\rho}_{ij}\}\right)$ to be the minimizer of (\ref{eq:lagrangian}) subject to the constraints that $\Tr[\rho_{i}]=1$ and $\tilde{\rho}_{ij}\succeq0$ for all $i<j$, holding $\{\rho_{ij},\Lambda_{ij},\Lambda_{ij}^{(1)},\Lambda_{ij}^{(2)}\}_{i<j}$ fixed} \label{eq:ADMM2}

\For{each pair $i<j$}
\State{$\Lambda_{ij}\leftarrow\Lambda_{ij}+\mu\left(\tilde{\rho}_{ij}-\rho_{ij}\right)$
}
\State{$\Lambda_{ij}^{(1)}\leftarrow\Lambda_{ij}^{(1)}+\nu(A_{1}[\rho_{ij}]-\rho_{i})$
}
\State{$\Lambda_{ij}^{(2)}\leftarrow\Lambda_{ij}^{(2)}+\nu(A_{2}[\rho_{ij}]-\rho_{j})$
}
\EndFor

\EndWhile
\end{algorithmic}

\label{alg:ADMM}
\end{algorithm}
Notice that in step \ref{eq:ADMM1} of Algorithm \ref{alg:ADMM} ,
the $\rho_{ij}$ are all determined \emph{independently} as the solutions
of decoupled optimization problems 
\begin{eqnarray*}
\rho_{ij} & \leftarrow & \underset{\rho_{ij}}{\mbox{argmin}}\Bigg\{\left\langle H_{ij}',\rho_{ij}\right\rangle _{\mathrm{F}}+\frac{\mu}{2}\Vert\rho_{ij}-\tilde{\rho}_{ij}\Vert_{\mathrm{F}}^{2}+\frac{\nu}{2}\Vert A_{1}[\rho_{ij}]-\rho_{i}\Vert_{\mathrm{F}}^{2}+\frac{\nu}{2}\Vert A_{2}[\rho_{ij}]-\rho_{j}\Vert_{\mathrm{F}}^{2}\\
 &  & \quad\quad\quad\quad\quad-\ \left\langle \Lambda_{ij},\rho_{ij}\right\rangle _{\mathrm{F}}+\left\langle \Lambda_{ij}^{(1)},A_{1}[\rho_{ij}]\right\rangle _{\mathrm{F}}+\left\langle \Lambda_{ij}^{(2)},A_{2}[\rho_{ij}]\right\rangle _{\mathrm{F}}\Bigg\}.
\end{eqnarray*}
 After suitable manipulation of the objective (neglecting constant
terms), we obtain 
\[
\frac{1}{2}\left\langle \rho_{ij},\left(\mu+\nu A_{1}^{*}A_{1}+\nu A_{2}^{*}A_{2}\right)\rho_{ij}\right\rangle _{\mathrm{F}}-\left\langle \mu\tilde{\rho}_{ij}+A_{1}^{*}\left[\nu\rho_{i}-\Lambda_{ij}^{(1)}\right]+A_{2}^{*}\left[\nu\rho_{j}-\Lambda_{ij}^{(2)}\right]+\Lambda_{ij}-H_{ij}',\rho_{ij}\right\rangle _{\mathrm{F}},
\]
 which can be exactly optimized via the update 
\begin{equation}
\rho_{ij}\leftarrow\left(\mu+\nu A_{1}^{*}A_{1}+\nu A_{2}^{*}A_{2}\right)^{-1}\left(\mu\tilde{\rho}_{ij}+A_{1}^{*}\left[\nu\rho_{i}-\Lambda_{ij}^{(1)}\right]+A_{2}^{*}\left[\nu\rho_{j}-\Lambda_{ij}^{(2)}\right]+\Lambda_{ij}-H_{ij}'\right).\label{eq:rhoijUpdate}
\end{equation}

Meanwhile, in step \ref{eq:ADMM2}, the $\rho_{i}$ and\emph{ }the
$\tilde{\rho}_{ij}$ can all be updated via decoupled optimization
problems. In particular, we find that 
\begin{equation}
\tilde{\rho}_{ij}\leftarrow\underset{\tilde{\rho}_{ij}\succeq0}{\mbox{argmin}}\left\{ \Vert\tilde{\rho}_{ij}-(\rho_{ij}-\mu^{-1}\Lambda_{ij})\Vert^{2}\right\} =\Pi_{\succeq0}\left(\rho_{ij}-\mu^{-1}\Lambda_{ij}\right).\label{eq:rhoijTildeUpdate}
\end{equation}
 Finally we turn to the $\rho_{i}$ update. Collecting the relevant
terms we have that 
\[
\rho_{i}\leftarrow\underset{\rho_{i}\,:\,\Tr[\rho_{i}]=1}{\mbox{argmin}}\left\{ \left\langle H_{i}',\rho_{i}\right\rangle _{\mathrm{F}}+\sum_{j>i}\left(\frac{\nu}{2}\Vert\rho_{i}-A_{1}[\rho_{ij}]\Vert_{\mathrm{F}}^{2}-\left\langle \Lambda_{ij}^{(1)},\rho_{i}\right\rangle _{\mathrm{F}}\right)+\sum_{j<i}\left(\frac{\nu}{2}\Vert\rho_{i}-A_{2}[\rho_{ji}]\Vert_{\mathrm{F}}^{2}-\left\langle \Lambda_{ji}^{(2)},\rho_{i}\right\rangle _{\mathrm{F}}\right)\right\} .
\]
 Observe that the objective may be rewritten as 
\[
\frac{(M-1)\nu}{2}\Vert\rho_{i}\Vert_{\mathrm{F}}^{2}-\left\langle \sum_{j>i}\left(\nu A_{1}[\rho_{ij}]+\Lambda_{ij}^{(1)}\right)+\sum_{j<i}\left(\nu A_{2}[\rho_{ji}]+\Lambda_{ji}^{(2)}\right)-H_{i}',\rho_{i}\right\rangle _{\mathrm{F}},
\]
 which we must minimize subject to $\Tr\left[\rho_{i}\right]=1$.
This is simply a constrained least squares problem, the solution of
which yields the update 
\begin{equation}
\rho_{i}\leftarrow\frac{1}{\nu(M-1)}\left(\sum_{j>i}\left(\nu A_{1}[\rho_{ij}]+\Lambda_{ij}^{(1)}\right)+\sum_{j<i}\left(\nu A_{2}[\rho_{ji}]+\Lambda_{ji}^{(2)}\right)-H_{i}'\right)+zI_{m},\label{eq:rhoiUpdate}
\end{equation}
 where $z$ is a Lagrange multiplier chosen to satisfy the constraint.

Then via (\ref{eq:rhoijUpdate}), (\ref{eq:rhoijTildeUpdate}), and
(\ref{eq:rhoiUpdate}), we can rewrite Algorithm \ref{alg:ADMM} concretely
as the equivalent Algorithm \ref{alg:ADMMconcrete}. Observe that
all of the for-loops in Algorithm \ref{alg:ADMMconcrete} can be run
in parallel. 
\begin{algorithm}
\caption{Details of the ADMM-type method for (\ref{eq:BP})}

\begin{algorithmic}[1]
\Require{$\mu,\nu > 0, \ \{H'_i, \rho_i\},\ \{H'_{ij}, \tilde{\rho}_{ij}, \Lambda_{ij}, \Lambda_{ij}^{(1)}, \Lambda_{ij}^{(2)}\}_{i<j}$}
\While{not converged}

\For{each pair $i<j$}
\State{$\rho_{ij}\leftarrow\left(\mu+\nu A_{1}^{*}A_{1}+\nu A_{2}^{*}A_{2}\right)^{-1}\left(\mu\tilde{\rho}_{ij}+A_{1}^{*}\left[\nu\rho_{i}-\Lambda_{ij}^{(1)}\right]+A_{2}^{*}\left[\nu\rho_{j}-\Lambda_{ij}^{(2)}\right]+\Lambda_{ij}-H_{ij}'\right)$
}
\EndFor

\For{each pair $i<j$}
\State{$\tilde{\rho}_{ij}\leftarrow \Pi_{\succeq0}\left(\rho_{ij}-\mu^{-1}\Lambda_{ij}\right)$
}
\EndFor
\For{each $i$}
\State{$\rho'_i \leftarrow \frac{1}{\nu(M-1)}\left(\sum_{j>i}\left(\nu A_{1}[\rho_{ij}]+\Lambda_{ij}^{(1)}\right)+\sum_{j<i}\left(\nu A_{2}[\rho_{ji}]+\Lambda_{ji}^{(2)}\right)-H_{i}'\right)$
}
\State{$z \leftarrow m^{-1}(1-\Tr[\rho'_i])$
}
\State{$\rho_i \leftarrow \rho'_i + z I_m$
}
\EndFor

\For{each pair $i<j$}
\State{$\Lambda_{ij}\leftarrow\Lambda_{ij}+\mu\left(\tilde{\rho}_{ij}-\rho_{ij}\right)$
}
\State{$\Lambda_{ij}^{(1)}\leftarrow\Lambda_{ij}^{(1)}+\nu(A_{1}[\rho_{ij}]-\rho_{i})$
}
\State{$\Lambda_{ij}^{(2)}\leftarrow\Lambda_{ij}^{(2)}+\nu(A_{2}[\rho_{ij}]-\rho_{j})$
}
\EndFor

\EndWhile
\end{algorithmic}

\label{alg:ADMMconcrete}
\end{algorithm}

\subsection{Practical scheme\label{sub:Practical-scheme}}

Now we return to the full Algorithm \ref{alg:gradAscent} where we optimize $X$ via an idealized gradient ascent. Instead of exactly implementing
step \ref{eq:exactBP} of Algorithm \ref{alg:gradAscent}, we replace
it with a single iteration of Algorithm \ref{alg:ADMMconcrete}, yielding
our practical approach Algorithm \ref{alg:practical} for solving
the two-marginal relaxation (\ref{eq:sdpObj}).

\begin{algorithm}
\caption{Practical ADMM / projected dual gradient ascent method for (\ref{eq:sdpObj})}

\begin{algorithmic}[1]
\Require{$\ve,\mu,\nu > 0, \ X\succeq 0,\  \{\rho_i\},\ \{\tilde{\rho}_{ij}, \Lambda_{ij}, \Lambda_{ij}^{(1)}, \Lambda_{ij}^{(2)}\}_{i<j}$}
\While{not converged}
\For{each $i$}
\State{$H'_{i} \leftarrow H_i[X_{ii}]$}
\EndFor
\For{each pair $i<j$}
\State{$H'_{ij} \leftarrow H_{ij}[X_{ij}]$}
\EndFor

\For{each pair $i<j$}
\State{$\rho_{ij}\leftarrow\left(\mu+\nu A_{1}^{*}A_{1}+\nu A_{2}^{*}A_{2}\right)^{-1}\left(\mu\tilde{\rho}_{ij}+A_{1}^{*}\left[\nu\rho_{i}-\Lambda_{ij}^{(1)}\right]+A_{2}^{*}\left[\nu\rho_{j}-\Lambda_{ij}^{(2)}\right]+\Lambda_{ij}-H_{ij}'\right)$
}
\EndFor

\For{each pair $i<j$}
\State{$\tilde{\rho}_{ij}\leftarrow \Pi_{\succeq0}\left(\rho_{ij}-\mu^{-1}\Lambda_{ij}\right)$
\label{eq:rhoijProj}}
\EndFor
\For{each $i$}
\State{$\rho'_i \leftarrow \frac{1}{\nu(M-1)}\left(\sum_{j>i}\left(\nu A_{1}[\rho_{ij}]+\Lambda_{ij}^{(1)}\right)+\sum_{j<i}\left(\nu A_{2}[\rho_{ji}]+\Lambda_{ji}^{(2)}\right)-H_{i}'\right)$
}
\State{$z \leftarrow m^{-1}(1-\Tr[\rho'_i])$
}
\State{$\rho_i \leftarrow \rho'_i + z I_m$
}
\EndFor

\For{each pair $i<j$}
\State{$\Lambda_{ij}\leftarrow\Lambda_{ij}+\mu\left(\tilde{\rho}_{ij}-\rho_{ij}\right)$
}
\State{$\Lambda_{ij}^{(1)}\leftarrow\Lambda_{ij}^{(1)}+\nu(A_{1}[\rho_{ij}]-\rho_{i})$
}
\State{$\Lambda_{ij}^{(2)}\leftarrow\Lambda_{ij}^{(2)}+\nu(A_{2}[\rho_{ij}]-\rho_{j})$
}
\EndFor

\State{$X\leftarrow\Pi_{\succeq0}\left(X+\ve G[\{\rho_{i}\},\{\rho_{ij}\}_{i<j}]\right)$} \label{eq:lastStep}

\EndWhile
\end{algorithmic}

\label{alg:practical}
\end{algorithm}

\subsection{Exploiting translation-invariance\label{sub:Exploiting-translation-invarianc}}

One of the most expensive step in Algorithm 5 is step 24, where a projection to the positive semidefinite cone is required. We now discuss how translation-invariance of a lattice system can
be exploited algorithmically in the solution of the two-marginal relaxation
(\ref{eq:sdpObj}). It is convenient in this section to adopt zero-indexing
for the site, i.e., to index the sites by the multi-index $\mathbf{i}=(i_{1},\ldots,i_{d})$
where $d$ is the lattice dimension and  $i_{p}=0,\ldots,M-1$.
Then we assume translation-invariance in that $\hat{H}_{\mathbf{i}}=\hat{H}_{\mathbf{0}}$
for all $i$ and $\hat{H}_{\mathbf{i}\mathbf{j}}=\hat{H}_{0,\mathbf{j}-\mathbf{i}}$
for all $i<j$. The symmetries $\rho_{\mathbf{i}}=\rho_{\mathbf{0}}$
and $\rho_{\mathbf{i}\mathbf{j}}=\rho_{0,\mathbf{j}-\mathbf{i}}$
are in turn guaranteed to be satisfied by some optimizer of (\ref{eq:sdpObj})~\cite{LinLindsey2020}.

Then we can implement Algorithm \ref{alg:practical} (whose iterations
preserve this symmetry) without any reference to variables besides
$\rho_{\mathbf{0}}$ and the $\rho_{\mathbf{0},\mathbf{j}}$. The
main challenge is the implementation of step \ref{eq:lastStep} of
Algorithm \ref{alg:practical}. Given $\rho_{\mathbf{0}}$ and the
$\rho_{\mathbf{0},\mathbf{j}}$, we can only compute the top row $X'_{\mathbf{0},\mathbf{j}}$
of $X':=X+\ve G[\{\rho_{i}\},\{\rho_{ij}\}_{i<j}]$. However, by translation
invariance, the rest of $X'$ is determined by the property that
$X'_{\mathbf{i},\mathbf{j}}=X'_{\mathbf{0},\mathbf{j}-\mathbf{i}}$. (In the case $d=1$, $X'$ is a \emph{block-circulant} matrix, though a more general term is lacking for the case of arbitrary $d$.)
Via translation invariance, $X'$ is block-diagonalized by
the block-discrete Fourier transform. More precisely, one can write $X'$ as $X' = (\mathcal{F}\otimes I) \hat{X}' (\mathcal{F}^* \otimes I)$, where $\mathcal{F}$ indicates the appropriate $d$-dimensional discrete Fourier transform matrix, $I$ the identity matrix of the appropriate block size, $\otimes$ the Kronecker product, and $\hat{X}'$ a block-diagonal matrix. Hence to compute the projection
$\Pi_{\succeq0}[X']$ we first compute the diagonal block $\hat{X}'_{\mathbf{k}}$ of $\hat{X}'$ as 
\[
\hat{X}'_{\mathbf{k}}=\frac{1}{\sqrt{M^{d}}}\sum_{j_{1},\ldots,j_{d}=0}^{M-1}\exp\left(-\iota\frac{2\pi\mathbf{j}\cdot\mathbf{k}}{M}\right)X'_{0,\mathbf{j}}.
\]
Note that the blocks $\hat{X}'_{\mathbf{k}}$ (concatenated into a block row) can be viewed as the entrywise discrete Fourier transform of the first block row of $X'$, hence can be computed simultaneously via FFT. Then project $Y_{\mathbf{k}}:=\Pi_{\succeq0}\left[\hat{X}'_{\mathbf{k}}\right]$ for all $\mathbf{k}$,
and set 
\[
X_{\mathbf{0},\mathbf{j}}\leftarrow\frac{1}{\sqrt{M^{d}}}\sum_{j_{1},\ldots,j_{d}=0}^{M-1}\exp\left(\iota\frac{2\pi\mathbf{j}\cdot\mathbf{k}}{M}\right)Y_{\mathbf{k}}.
\]
 The final pseudocode for the translation-invariant setting is given
in Algorithm \ref{alg:practicalTI}.

\begin{algorithm}
\caption{Translation-invariant ADMM / projected dual gradient ascent method
for (\ref{eq:sdpObj})}

\begin{algorithmic}[1]
\Require{$\ve,\mu,\nu > 0, \ \left(X_{\mathbf{0},\mathbf{j}} \right),\  \rho_{\mathbf{0}},\ \{\tilde{\rho}_{\mathbf{0},\mathbf{j}}, \Lambda_{\mathbf{0},\mathbf{j}}, \Lambda_{\mathbf{0},\mathbf{j}}^{(1)}, \Lambda_{\mathbf{0},\mathbf{j}}^{(2)}\}_{\mathbf{j} \neq \mathbf{0}}$}
\While{not converged}

\State{$H'_{\mathbf{0}} \leftarrow H_{
\mathbf{0}}[X_{\mathbf{0},\mathbf{0}}]$}

\For{each $\mathbf{j} \neq \mathbf{0}$}
\State{$H'_{\mathbf{0},\mathbf{j}} \leftarrow H_{\mathbf{0},\mathbf{j}}[X_{\mathbf{0},\mathbf{j}}]$}
\EndFor

\For{each $\mathbf{j} \neq \mathbf{0}$}
\State{$\rho_{\mathbf{0},\mathbf{j}}\leftarrow\left(\mu+\nu A_{1}^{*}A_{1}+\nu A_{2}^{*}A_{2}\right)^{-1}\left(\mu\tilde{\rho}_{\mathbf{0},\mathbf{j}}+A_{1}^{*}\left[\nu\rho_{\mathbf{0}}-\Lambda_{\mathbf{0},\mathbf{j}}^{(1)}\right]+A_{2}^{*}\left[\nu\rho_{\mathbf{0}}-\Lambda_{\mathbf{0},\mathbf{j}}^{(2)}\right]+\Lambda_{\mathbf{0},\mathbf{j}}-H_{\mathbf{0},\mathbf{j}}'\right)$
}
\EndFor

\For{each $\mathbf{j} \neq \mathbf{0}$}
\State{$\tilde{\rho}_{\mathbf{0},\mathbf{j}}\leftarrow \Pi_{\succeq0}\left(\rho_{\mathbf{0},\mathbf{j}}-\mu^{-1}\Lambda_{\mathbf{0},\mathbf{j}}\right)$
}
\EndFor

\State{$\rho'_{\mathbf{0}} \leftarrow \frac{1}{\nu(M-1)}\left(\sum_{\mathbf{j}\neq \mathbf{0}}\left(\nu A_{1}[\rho_{\mathbf{0},\mathbf{j}}]+\Lambda_{\mathbf{0},\mathbf{j}}^{(1)}\right)-H_{\mathbf{0}}'\right)$
}
\State{$z \leftarrow m^{-1}(1-\Tr[\rho'_{\mathbf{0}}])$
}
\State{$\rho_{\mathbf{0}} \leftarrow \rho'_{\mathbf{0}} + z I_m$
}

\For{each $\mathbf{j} \neq \mathbf{0}$}
\State{$\Lambda_{\mathbf{0},\mathbf{j}}\leftarrow\Lambda_{\mathbf{0},\mathbf{j}}+\mu\left(\tilde{\rho}_{\mathbf{0},\mathbf{j}}-\rho_{\mathbf{0},\mathbf{j}}\right)$
}
\State{$\Lambda_{\mathbf{0},\mathbf{j}}^{(1)}\leftarrow\Lambda_{\mathbf{0},\mathbf{j}}^{(1)}+\nu(A_{1}[\rho_{\mathbf{0},\mathbf{j}}]-\rho_{\mathbf{0}})$
}
\State{$\Lambda_{\mathbf{0},\mathbf{j}}^{(2)}\leftarrow\Lambda_{\mathbf{0},\mathbf{j}}^{(2)}+\nu(A_{2}[\rho_{\mathbf{0},\mathbf{j}}]-\rho_{\mathbf{0}})$
}
\EndFor

\State{$X_{\mathbf{0},\mathbf{0}}'\leftarrow X_{\mathbf{0},\mathbf{0}}+\ve G_{\mathbf{0},\mathbf{0}}\left[\rho_{\mathbf{0}}\right]$ \label{eq:Xbegin}
}
\For{each $\mathbf{j} \neq \mathbf{0}$}
\State{$X_{\mathbf{0},\mathbf{j}}'\leftarrow X_{\mathbf{0},\mathbf{j}}+\ve G_{\mathbf{0},\mathbf{j}}\left[\rho_{\mathbf{0},\mathbf{j}}\right]$
}
\EndFor
\For{each $\mathbf{k}$}
\State{$\hat{X}'_{\mathbf{k}}=\frac{1}{\sqrt{M^{d}}}\sum_{\mathbf{j}}\exp\left(-\iota\frac{2\pi\mathbf{j}\cdot\mathbf{k}}{M}\right)X'_{0,\mathbf{j}}$
}
\State{$Y_{\mathbf{k}}\leftarrow \Pi_{\succeq0}\left[\hat{X}'_{\mathbf{k}}\right]$ \label{eq:Yproj}
}
\EndFor
\For{each $\mathbf{j}$}
\State{$X_{\mathbf{0},\mathbf{j}}\leftarrow\frac{1}{\sqrt{M^{d}}}\sum_{\mathbf{j}}\exp\left(\iota\frac{2\pi\mathbf{j}\cdot\mathbf{k}}{M}\right)Y_{\mathbf{k}}$ \label{eq:Xend}
}
\EndFor

\EndWhile
\end{algorithmic}

\label{alg:practicalTI}
\end{algorithm}

\subsection{Discussion of scaling\label{sub:Discussion-of-scaling}}

In Algorithm \ref{alg:practical}, observe that the for-loops run
over $M(M-1)$ pairs $i<j$, and the scaling bottleneck among these
loops is the projection $\tilde{\rho}_{ij}\leftarrow\Pi_{\succeq0}\left(\rho_{ij}-\mu^{-1}\Lambda_{ij}\right)$
occuring in step \ref{eq:rhoijProj}. Since this step requires full
diagonalization of a matrix of size $m^{2}\times m^{2}$, for which
the cost is $O(m^{6})$. Meanwhile, suppose for simplicity that $n_{i}=m^{2}$,
corresponding to the complete choice of operator collection $\left\{ O_{\alpha,i}\,:\,\alpha=n_{i}\right\} $
for each site. Then the size of $G$ and $X$ is $Mm^{2}\times Mm^{2}$.
Hence step \ref{eq:lastStep}, which involves a complete diagonalization
of a matrix of this size, costs $O(m^{6}M^{3})$, dominating the $O(m^{6}M^{2})$
cost of the for-loops. If our sites are in fact supersites, each formed
from clusters of $L$ sites in an underlying spin-$\frac{1}{2}$ model,
then $m=2^{L}$. Therefore the scaling is $O(2^{6L}M^{3})$ per iteration,
where $L$ is the cluster size and $M$ is the number of clusters.

Meanwhile, in the translation-invariant setting of Algorithm \ref{alg:practicalTI},
the for-loops run only over $O(M)$ sites, so---neglecting the update
for $X$---the asymptotic cost per iteration is $O(m^{6}M)$. Meanwhile,
the construction of the $\hat{X}_{\mathbf{k}}'$ in terms of the $X_{0,\mathbf{j}}$
can be achieved in time $O(m^{4}M\log M)$ via FFT. (Note that we
simply treat the lattice dimension $d$ as constant.) The cost of
each projection of step \ref{eq:Yproj} is $O(m^{6})$ via full diagonalization,
and forming the $X_{\mathbf{0},\mathbf{j}}$ in terms of the $Y_{\mathbf{k}}$
also costs $O(dm^{4}M\log M)$ in total via FFT. Hence the cost of
updating $X$ (i.e., steps \ref{eq:Xbegin} through \ref{eq:Xend})
is $O(m^{4}M\log M+m^{6}M)$. Hence the total cost per iteration of
Algorithm \ref{alg:practicalTI}, under the assumption that the sites
are supersites each composed of $L$ spin-$\frac{1}{2}$ sites, is
$O(2^{4L}M\log M+2^{6L}M)$, where $L$ is the cluster size and $M$
is the number of clusters.

The exponential scaling in the cluster size is unavoidable in our
formulation due to the exact treatment of reduced density operators
on the clusters (i.e., the cluster marginals). In this work we consider
clusters of size no larger than size $L=4$. Future work will investigate
the possibility of treating larger clusters by introducing further
relaxation and/or compression of the optimization variables to avoid
exponential scaling in the cluster size.

As a final comment, observe that every for-loop in Algorithms \ref{alg:practical}
and \ref{alg:practicalTI} can be run fully in parallel.

\section{Numerical experiments\label{sec:Numerical-experiments}}

The numerical experiments were implemented in $\textsf{MATLAB}$ following
Algorithm \ref{alg:practicalTI}. (We shall consider only translation-invariant
Hamiltonians.) We present results for the transverse-field Ising (TFI)
model (\ref{eq:TFI}), the anti-ferromagnetic Heisenberg (AFH) model
(\ref{eq:AFH}), the spinless fermion (SF) model (\ref{eq:spinlessFerm}),
and the long-range spinless fermion (LRSF) model (\ref{eq:spinlessFermLR}).
Note that Algorithm \ref{alg:practicalTI} can be applied in the fermionic
case \emph{mutatis mutandi} to the problem (\ref{eq:sdpFerm}). Throughout
we fix the value of the algorithmic parameters to be $\mu=\nu=10$,
$\ve=2$ (i.e., we do not tune them specifically to different problems).
The dual variables $\{ \Lambda_{ij}, \Lambda_{ij}^{(1)}, \Lambda_{ij}^{(2)}\}_{i\neq j}$ are all initialized to be zero, and the primal
density operator variables $\{ \rho_{i} \},\ \{\tilde{\rho}_{ij} \}_{i\neq j}$ are all initialized as multiples of the
identity with unit trace. $X$ is initialized as the identity. We
run Algorithm \ref{alg:practicalTI} for 10,000 iterations. (The convergence
behavior will be studied in detail below.)

\begin{figure}
\noindent \centering{}\includegraphics[bb=40bp 0bp 520bp 420bp,scale=0.45]{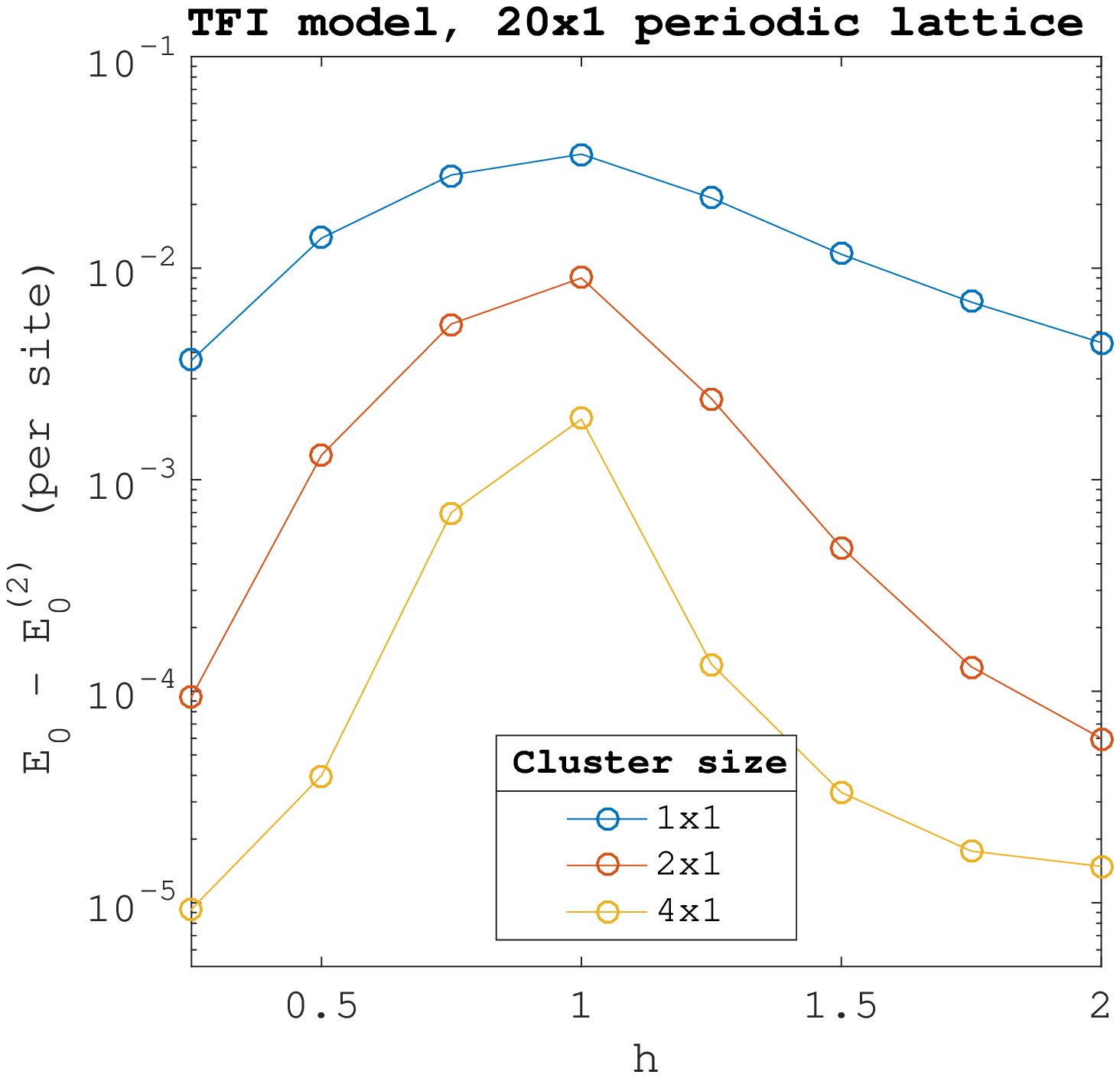}\includegraphics[bb=40bp 0bp 520bp 420bp,scale=0.45]{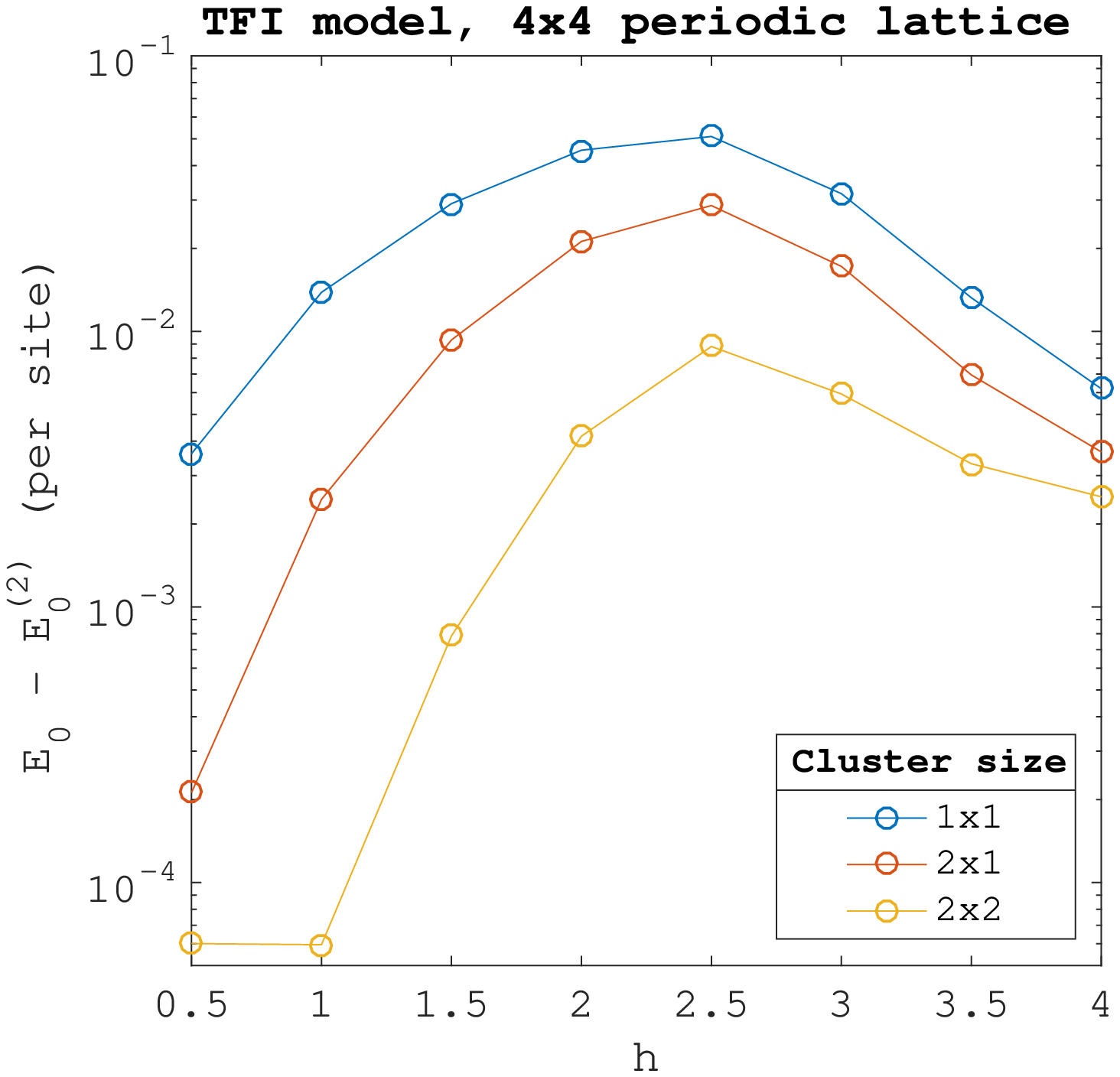}\caption{Relaxation error per site for the TFI model on $20\times1$ and $4\times4$
periodic lattices, pictured left and right respectively. Several cluster
sizes are considered according to the legends. Note that the relaxation
is exact at $h=0$ (not pictured).}
\label{fig:tfi_benchmark}
\end{figure}

First we consider the TFI model on a periodic $20\times1$ lattice,
which is small enough to be solved by exact diagonalization of (\ref{eq:TFI}).
We benchmark the per-site energy error of the two-marginal relaxation
with clusters of size $1\times1$, $2\times1$, and $4\times1$. (In these cases the semidefinite matrix variables $\rho_{ij}$ are \emph{each} of size 4, 16, and 256, respectively; refer to Section \ref{sub:Discussion-of-scaling} for further discussion of scaling.). The
results are shown in Figure \ref{fig:tfi_benchmark}. Observe that
the approximations yield lower bounds for the energy as the theory
requires, and these lower bounds become tighter as the cluster size
is increased. In the same figure we also consider the TFI model on
a periodic $4\times4$ lattice. Here we benchmark the energy error
of the two-marginal relaxation with clusters of size $1\times1$,
$2\times1$, and $2\times2$.

We perform completely analogous experiments for the AFH model with
similar conclusions. The results are shown in Tables \ref{tab:afh_benchmark_20x1}
and \ref{tab:afh_benchmark_4x4}.

\begin{table}

\centering{}%
\begin{tabular}{|c|c|c|}
\hline 
$1\times1$ clusters & $2\times1$ clusters & $4\times1$ clusters\tabularnewline
\hline 
0.5383 & 0.0521 & 0.0034\tabularnewline
\hline 
\end{tabular} \caption{Relaxation error per site for the AFH model on a $20\times1$ periodic
lattice for various cluster sizes. }
\label{tab:afh_benchmark_20x1}
\end{table}
\begin{table}
\centering{}%
\begin{tabular}{|c|c|c|}
\hline 
$1\times1$ clusters & $2\times1$ clusters & $2\times2$ clusters\tabularnewline
\hline 
0.6634 & 0.1851 & 0.0034\tabularnewline
\hline 
\end{tabular} \caption{Relaxation error per site for the AFH model on a $4\times4$ periodic
lattice for various cluster sizes. }
\label{tab:afh_benchmark_4x4}
\end{table}

\begin{figure}
\noindent \centering{}\includegraphics[bb=40bp 0bp 520bp 420bp,scale=0.45]{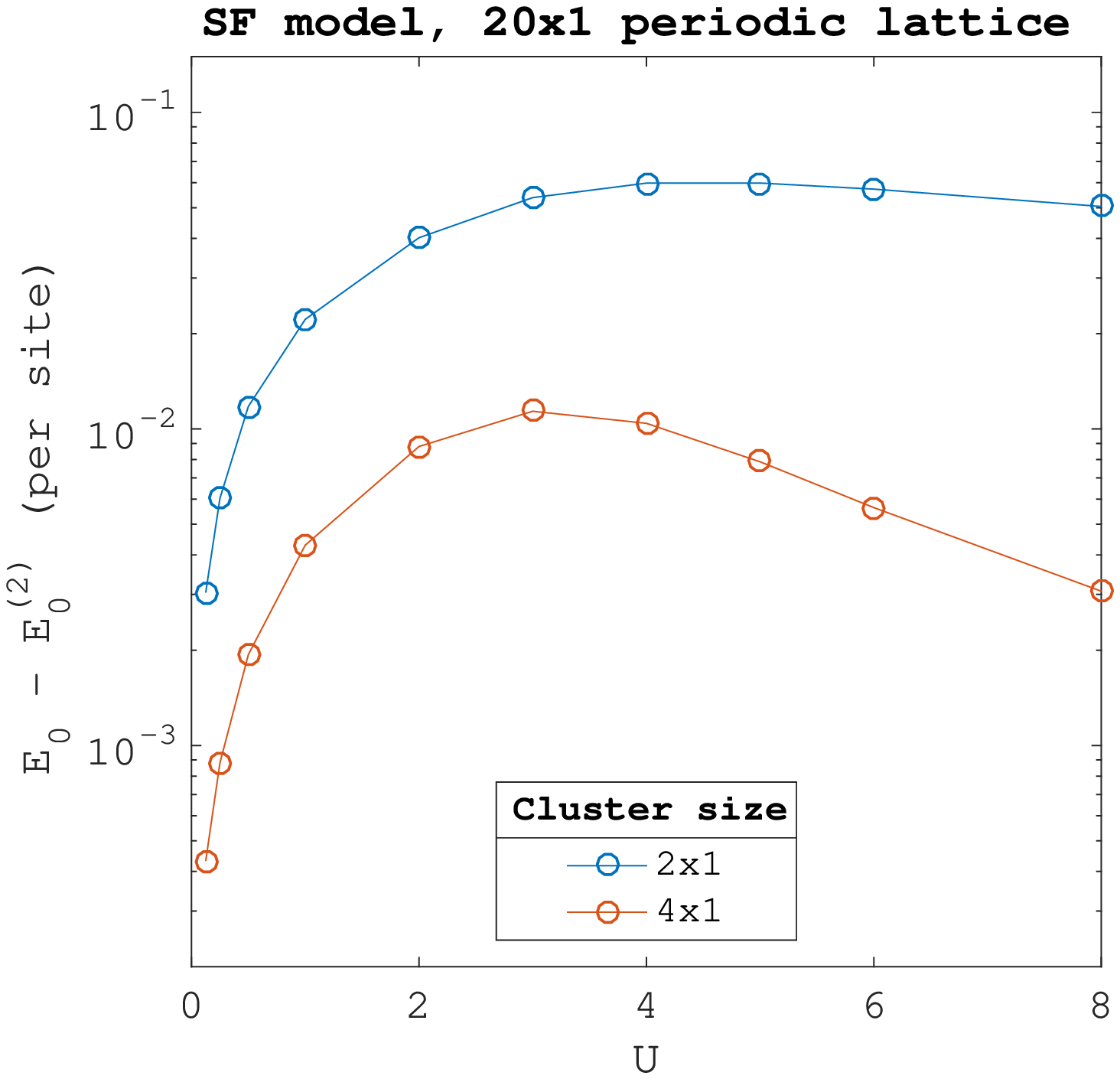}\includegraphics[bb=40bp 0bp 520bp 420bp,scale=0.45]{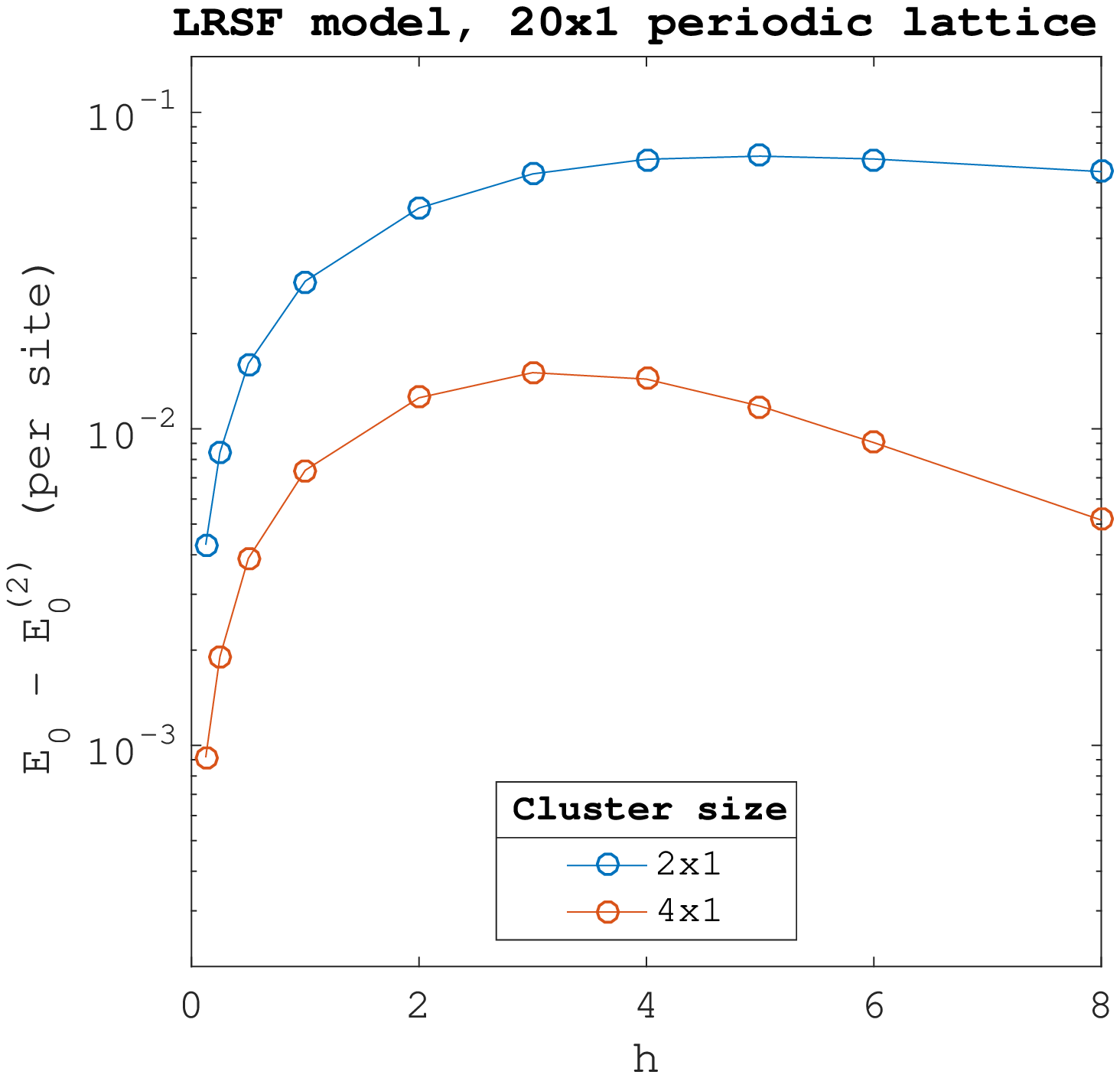}\caption{Relaxation error for the (short-range) spinless fermion \eqref{eq:spinlessFerm} and long-range
spinless fermion \eqref{eq:spinlessFermLR} models on a $20\times1$ periodic lattice, pictured
left and right respectively. Several cluster sizes are considered
according to the legends. Note that the relaxation is exact at $U=0$
(not pictured).}
\label{fig:sf_benchmark}
\end{figure}

We also benchmark the SF and LRSF models on a $20\times1$ periodic
lattice, with results pictured in Figure \ref{fig:sf_benchmark}.
Note that the fermionic relaxations are exact for $U=0$, as guaranteed
in~\cite{LinLindsey2020}.

Next we consider the TFI model on a periodic 100 \texttimes{} 1 lattice
for $h=0.5,1,1.5$. This problem is too large to solve by exact diagonalization.
We report the relaxation energy for several cluster sizes in Table
\ref{tab:tfi_100x1_energies}.

\begin{table}
\centering{}%
\begin{tabular}{|c|c|c|c|}
\hline 
 & $1\times1$ clusters & $2\times1$ clusters & $4\times1$ clusters\tabularnewline
\hline 
$h=0.5$ & $-1.0763$ & $-1.0648$ & $-1.0636$\tabularnewline
\hline 
$h=1$ & $-1.3084$ & $-1.2829$ & $-1.2761$\tabularnewline
\hline 
$h=1.5$ & $-1.6835$ & $-1.6724$ & $-1.6720$\tabularnewline
\hline 
\end{tabular} \caption{Relaxation energy per site for the TFI model with $h=0.5,1,1.5$ on
a $100\times1$ periodic lattice for various cluster sizes. }
\label{tab:tfi_100x1_energies}
\end{table}

\begin{figure}
\noindent \begin{centering}
\includegraphics[bb=40bp 0bp 500bp 420bp,clip,scale=0.33]{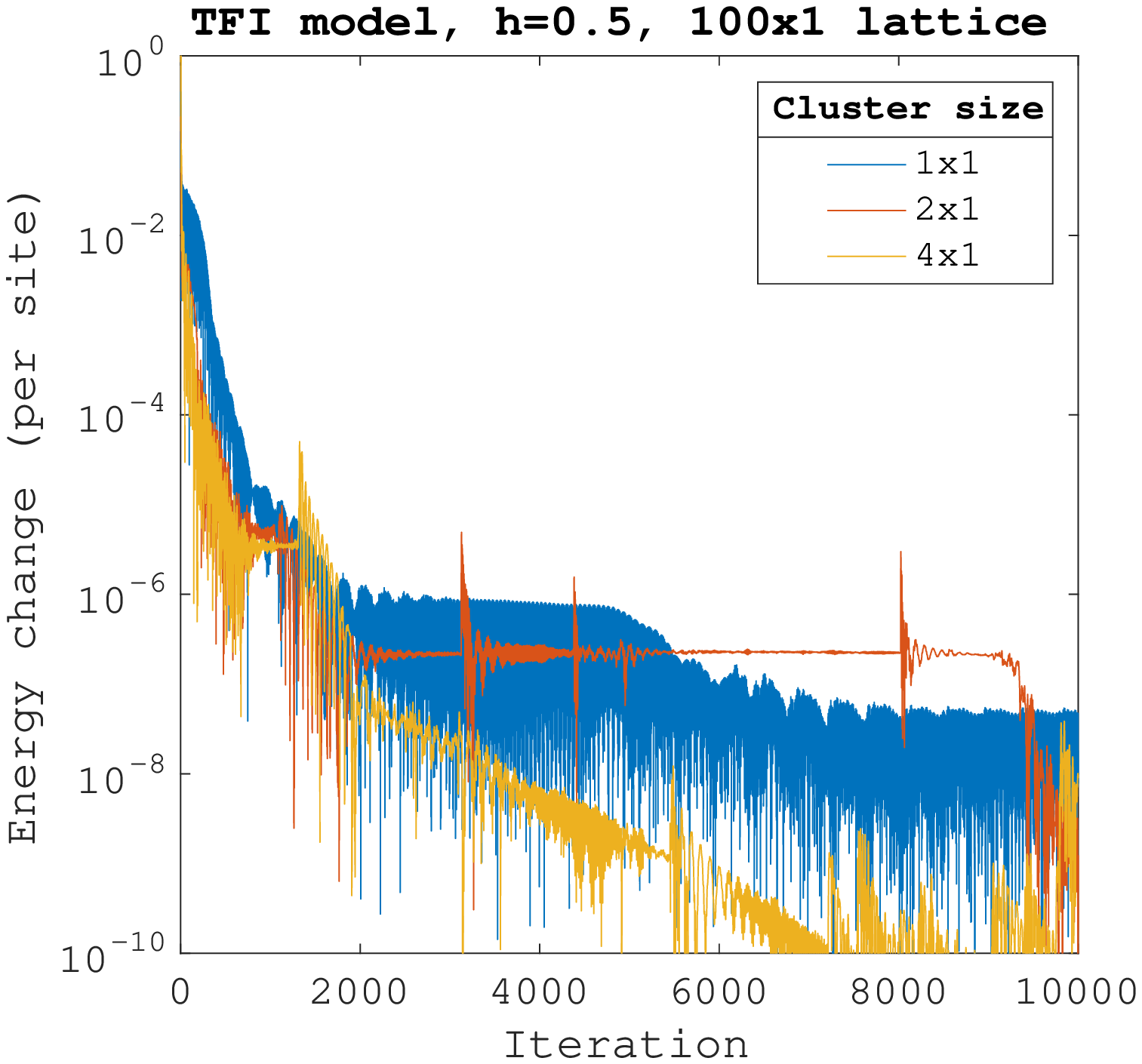}\includegraphics[bb=75bp 0bp 500bp 420bp,clip,scale=0.33]{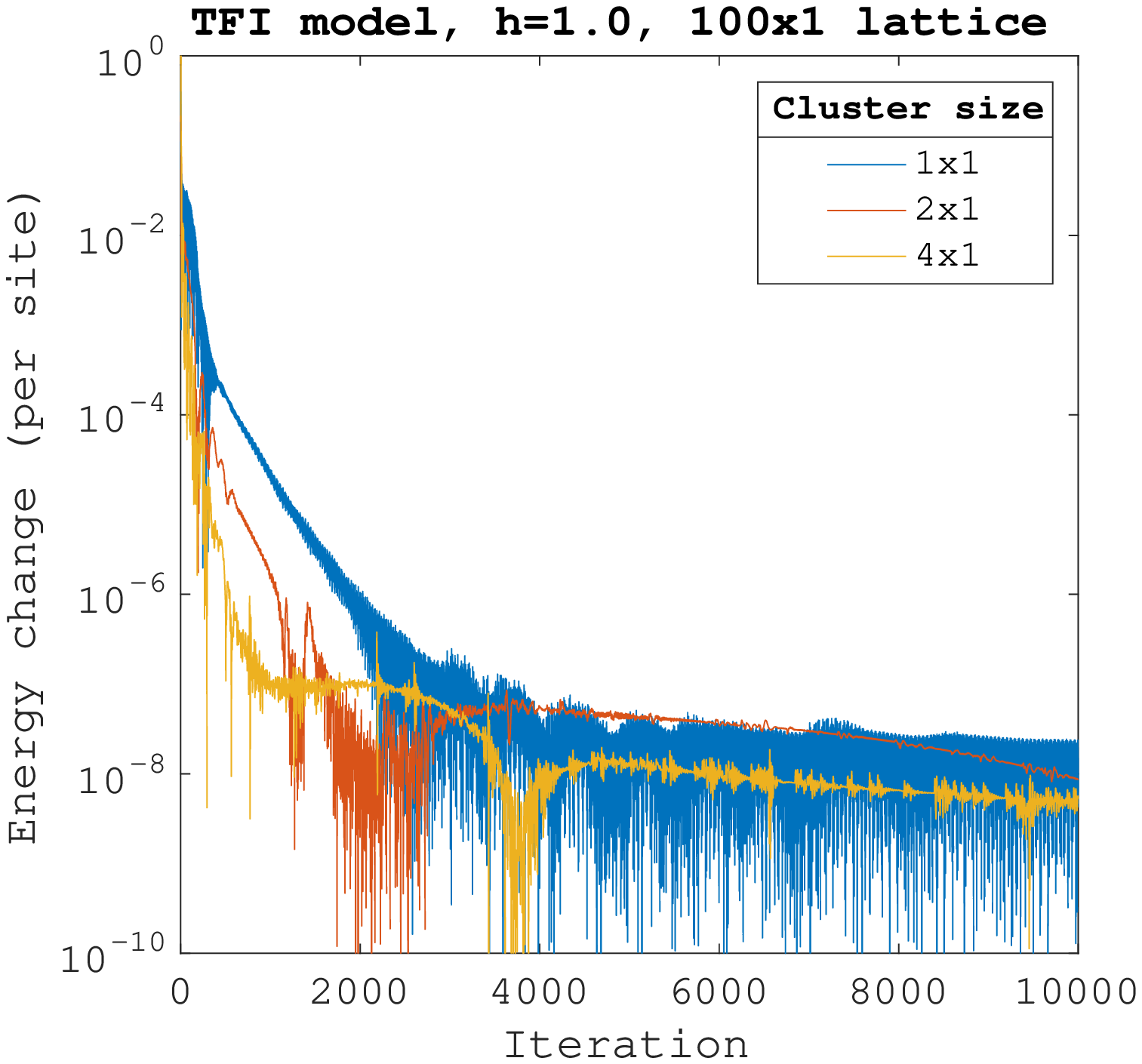}\includegraphics[bb=75bp 0bp 500bp 420bp,clip,scale=0.33]{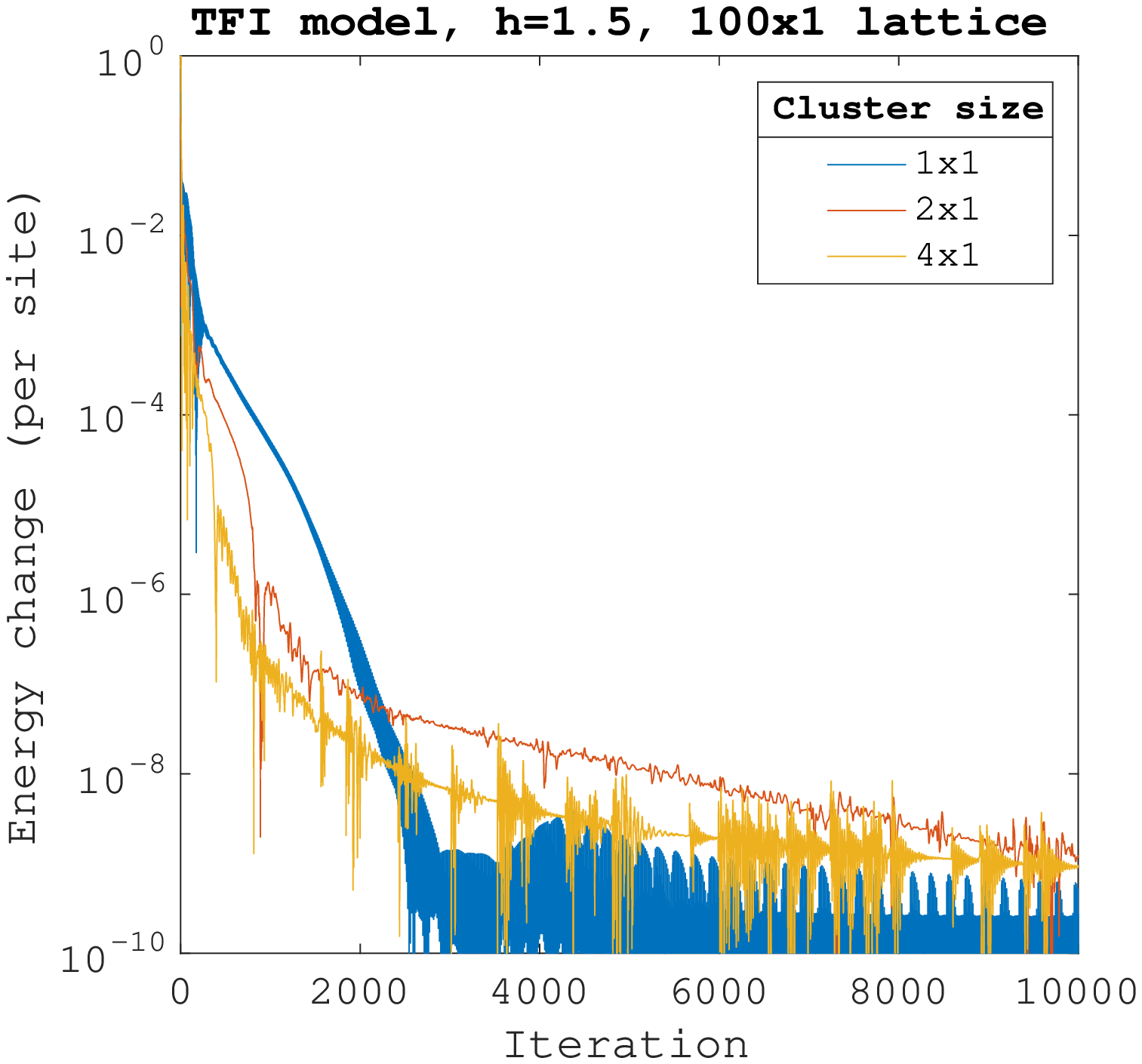}\vspace{2mm}
\par\end{centering}

\noindent \centering{}\includegraphics[bb=40bp 0bp 500bp 420bp,clip,scale=0.33]{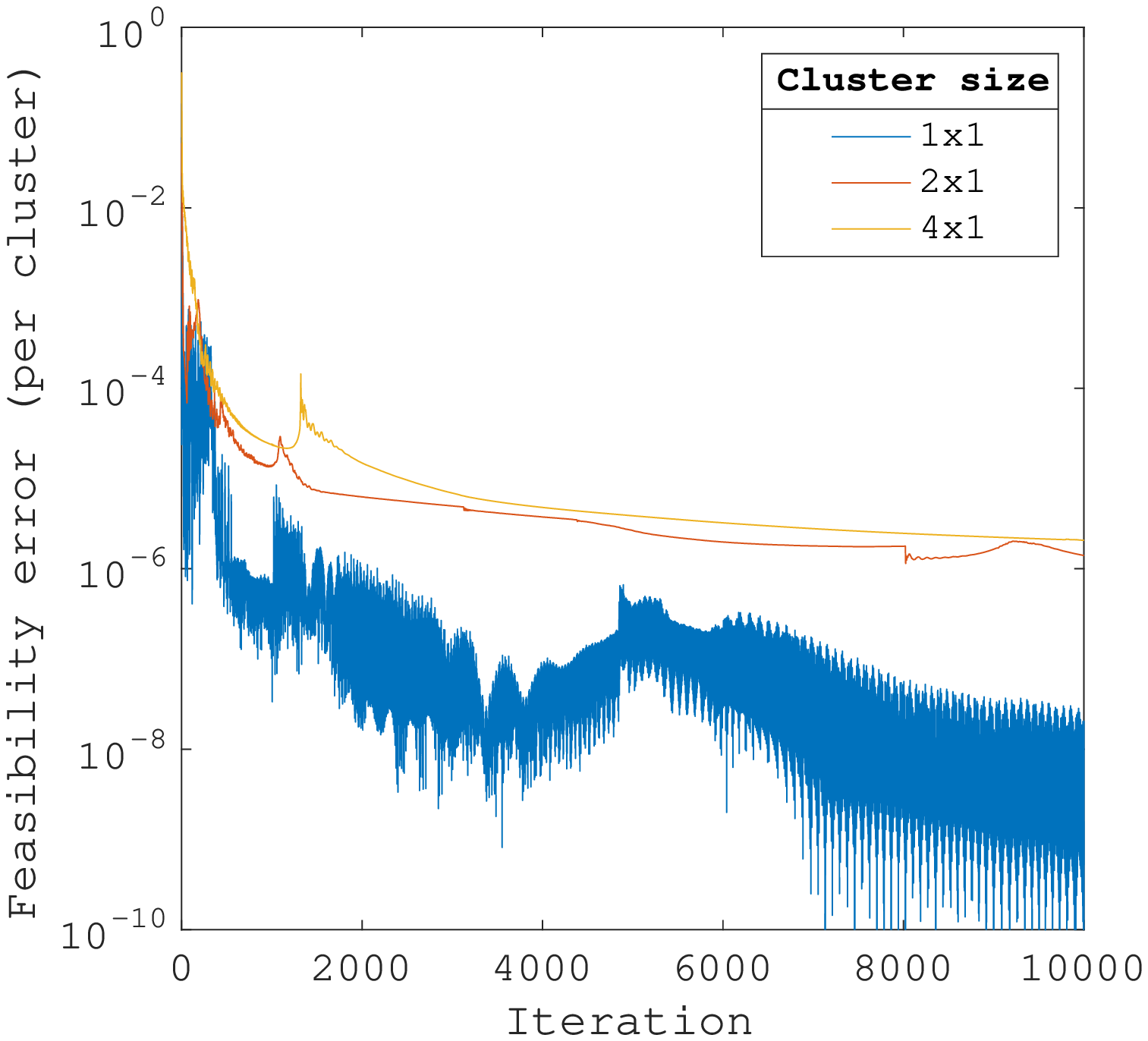}\includegraphics[bb=75bp 0bp 500bp 420bp,clip,scale=0.33]{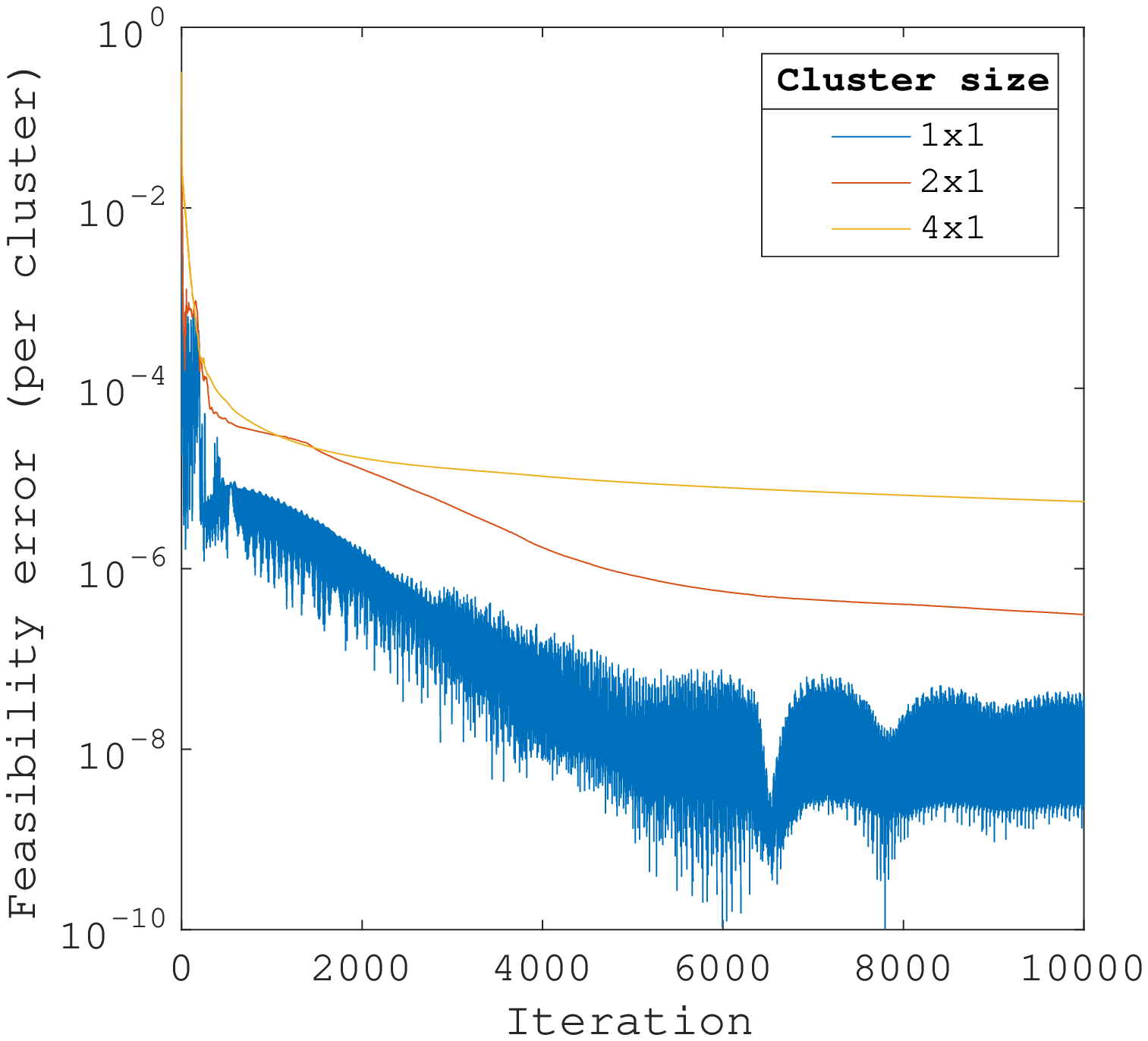}\includegraphics[bb=75bp 0bp 500bp 420bp,clip,scale=0.33]{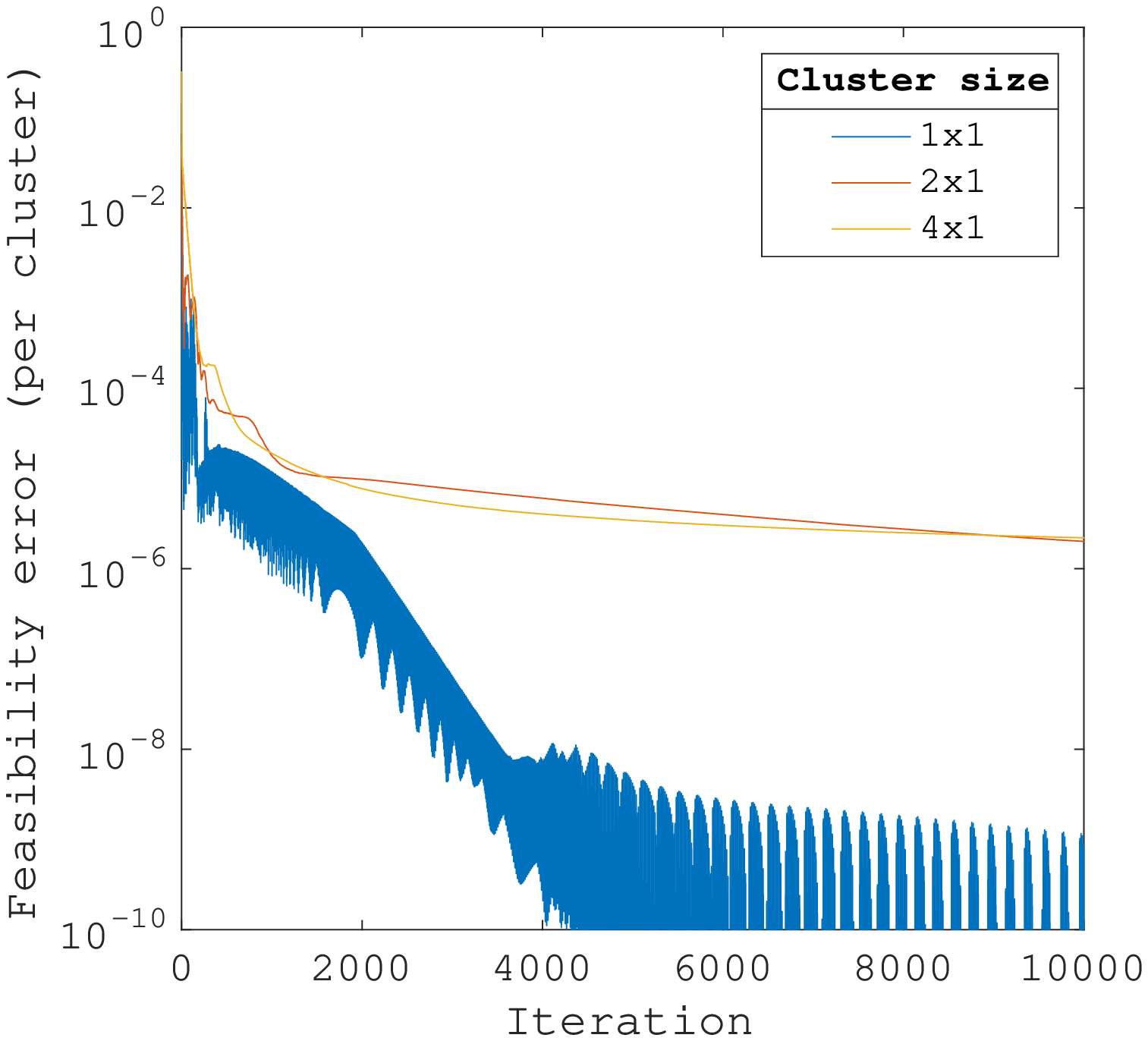}\caption{Per-site iteration-over-iteration energy change for $100\times1$
periodic TFI model with $h=0.5,1,1.5$ and several cluster sizes (top
row). Per-cluster feasibility error (\ref{eq:feasErr}) for same problems
(bottom row).}
\label{fig:tfi_100x1_conv_cluster}
\end{figure}

We track convergence behavior of Algorithm \ref{alg:practicalTI}
on these same problems. We use two different quantities to track convergence.
The first is the per-site primal objective (i.e., energy) change between
subsequent iterations. The second is the per-cluster feasibility error
for the equality constraints, defined as 
\begin{equation}
\sqrt{\frac{1}{M-1}\sum_{\mathbf{j}\neq0}\left(\Vert A_{1}[\rho_{\mathbf{0},\mathbf{j}}]-\rho_{\mathbf{0}}\Vert_{\mathrm{F}}^{2}+\Vert A_{2}[\rho_{\mathbf{0},\mathbf{j}}]-\rho_{\mathbf{0}}\Vert_{\mathrm{F}}^{2}+\Vert\rho_{\mathbf{0},\mathbf{j}}-\tilde{\rho}_{\mathbf{0},\mathbf{j}}\Vert_{\mathrm{F}}^{2}\right)}.\label{eq:feasErr}
\end{equation}
 We plot these quantities as functions of the iteration count in Figure
\ref{fig:tfi_100x1_conv_cluster}. It is possible that tuning the
parameters $\mu,\nu,\ve$ to a specific problem and specific choice
of clusters could yield smoother convergence profiles. However, even
using our fixed choice for all problems, we achieve convergence of
the per-site energy within $10^{-6}$ (which is dominated by the relaxation
error itself) in a number of iterations that does not seem to grow
with the cluster size.

\begin{figure}
\noindent \begin{centering}
\includegraphics[bb=40bp 0bp 500bp 420bp,clip,scale=0.33]{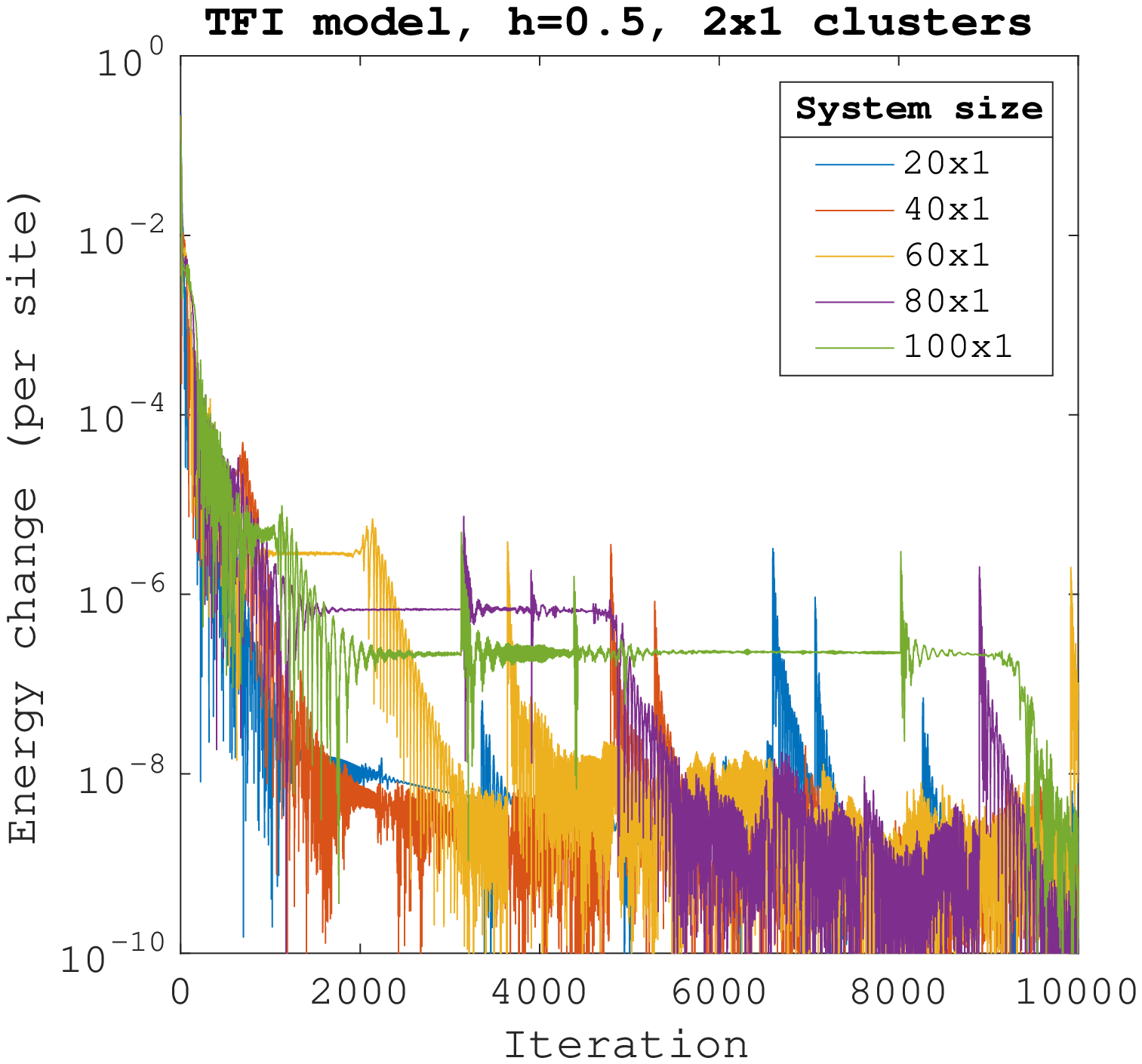}\includegraphics[bb=75bp 0bp 500bp 420bp,clip,scale=0.33]{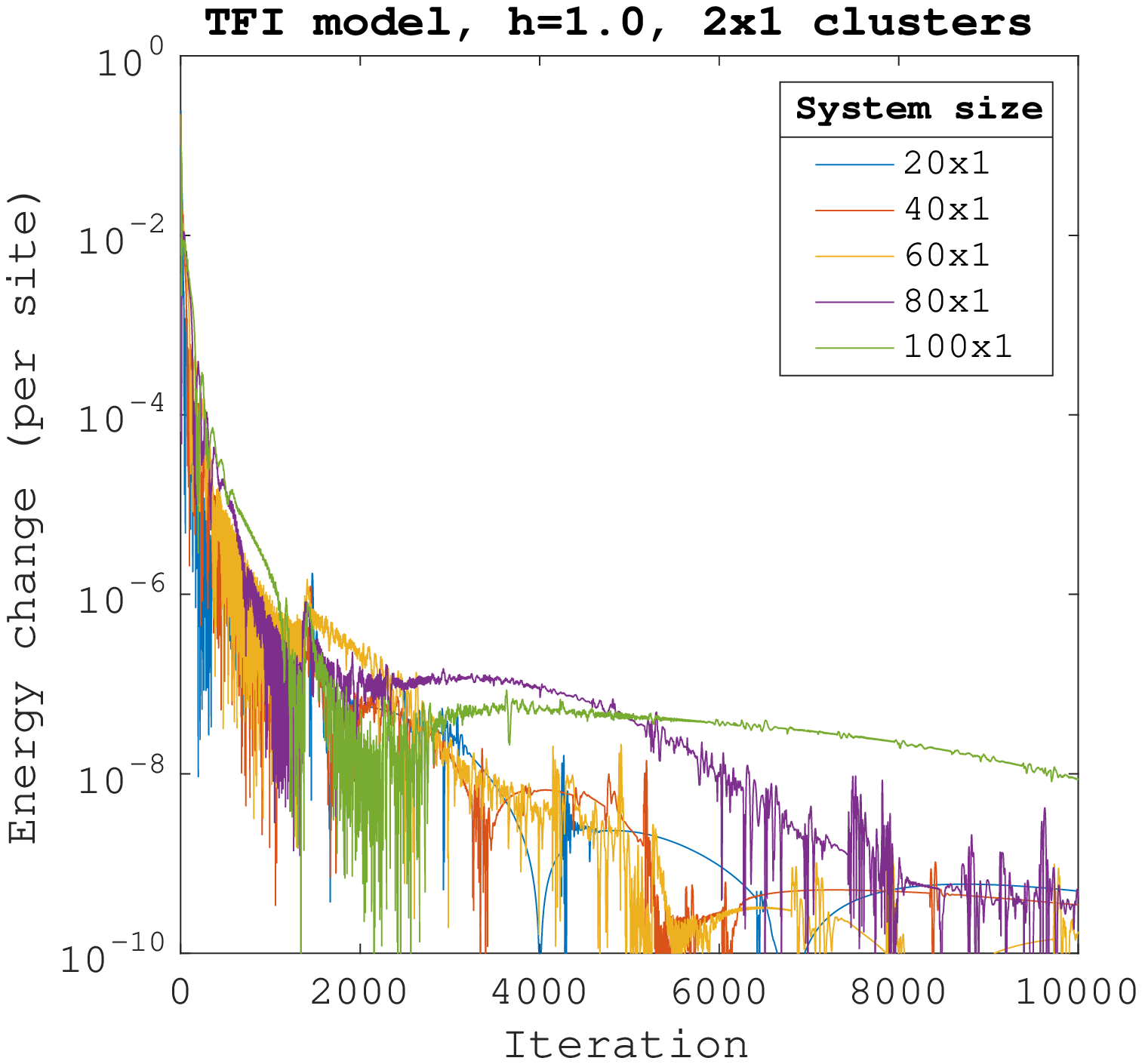}\includegraphics[bb=75bp 0bp 500bp 420bp,clip,scale=0.33]{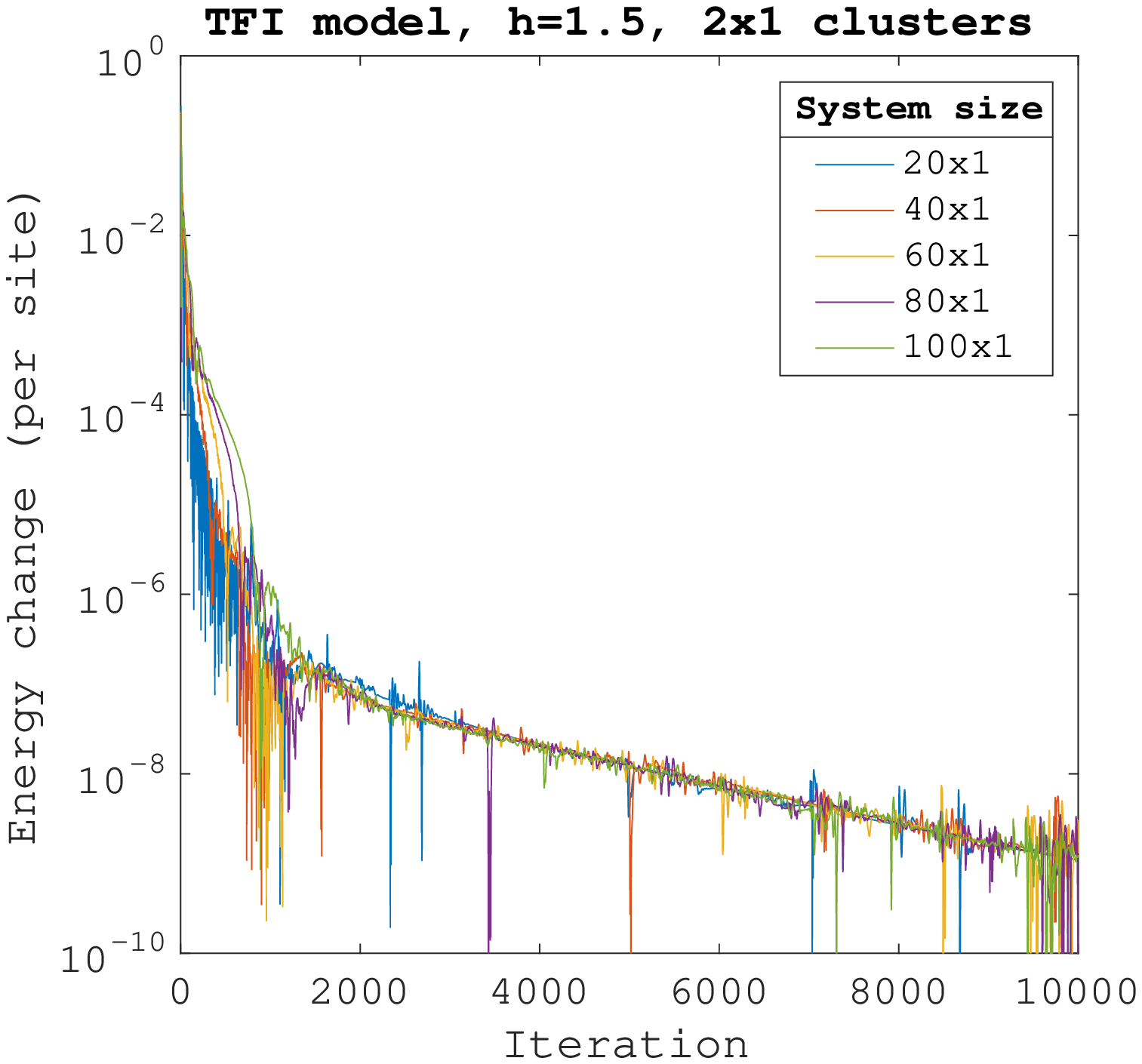}\vspace{2mm}
\par\end{centering}

\noindent \centering{}\includegraphics[bb=40bp 0bp 500bp 420bp,clip,scale=0.33]{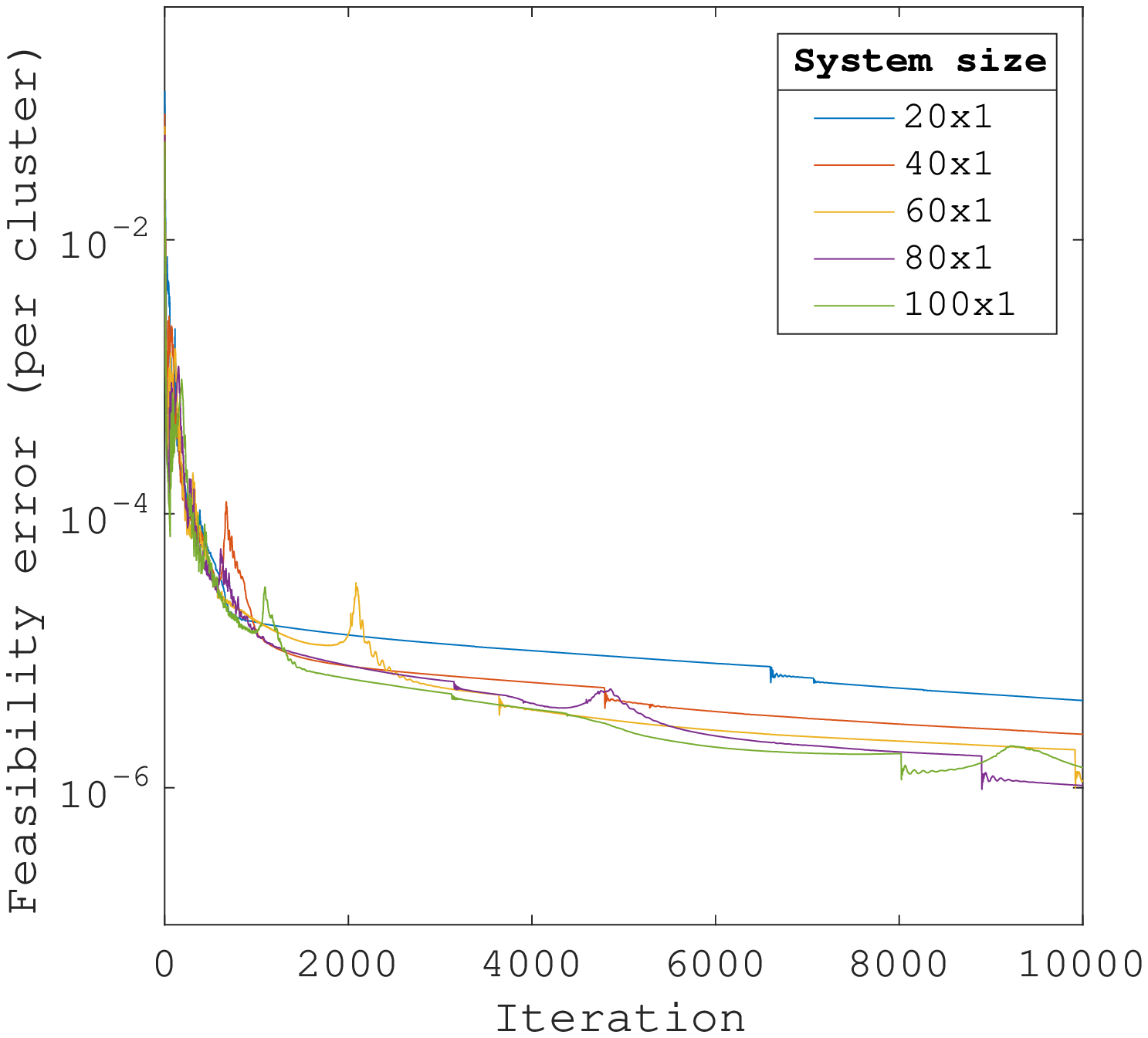}\includegraphics[bb=75bp 0bp 500bp 420bp,clip,scale=0.33]{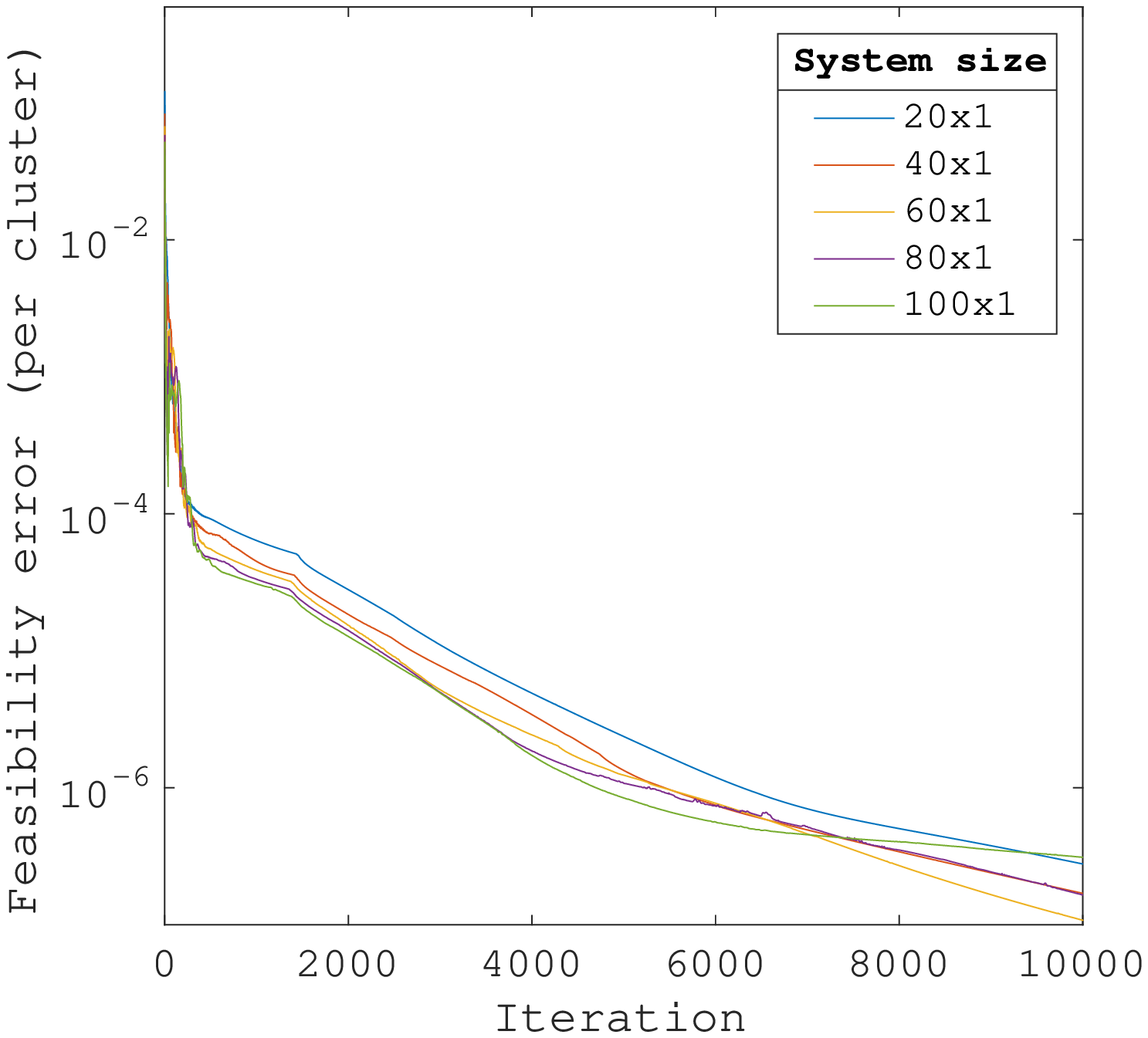}\includegraphics[bb=75bp 0bp 500bp 420bp,clip,scale=0.33]{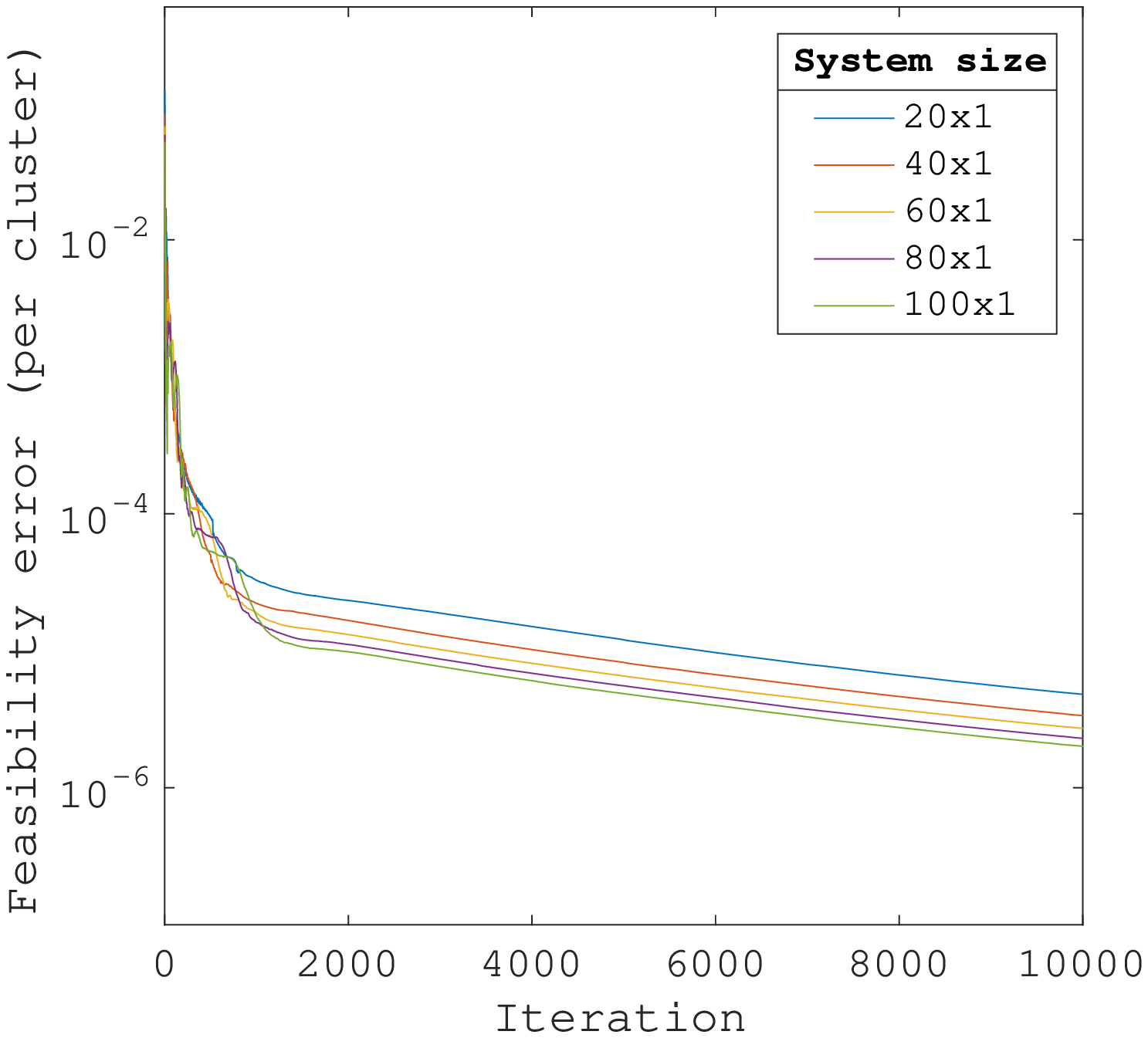}\caption{Per-site iteration-over-iteration energy change for periodic TFI models
of different sizes with $h=0.5,1,1.5$ and fixed $2\times1$ cluster
size (top row). Per-cluster feasibility error (\ref{eq:feasErr})
for same problems (bottom row).}
\label{fig:tfi_100x1_conv_size}
\end{figure}

Then we fix clusters of size $2\times1$ and vary the system size
of the TFI model with $h=0.5,1,1.5$, to investigate the effect of
system size on convergence. The results are shown in Figure \ref{fig:tfi_100x1_conv_size}.
The system size does not appear to have any obvious detrimental effect
on the convergence rate. 

To conclude, we report the relaxation energies obtained from these
last experiments in Table \ref{tab:tfi_size_energies}. We observe
that the relaxation energy approaches a limiting value in this thermodynamic limit, i.e., limit of infinite volume.

\begin{table}
\centering{}%
\begin{tabular}{|c|c|c|c|c|c|}
\hline 
 & $20\times1$ lattice & $40\times1$ lattice & $60\times1$ lattice & $80\times1$ lattice & $100\times1$ lattice\tabularnewline
\hline 
$h=0.5$ & $-1.064851$ & $-1.064795$ & $-1.064786$ & $-1.064779$ & $-1.064776$\tabularnewline
\hline 
$h=1$ & $-1.283534$ & $-1.283083$ & $-1.283003$ & $-1.282975$ & $-1.282949$\tabularnewline
\hline 
$h=1.5$ & $-1.672407$ & $-1.672394$ & $-1.672393$ & $-1.672393$ & $-1.672394$\tabularnewline
\hline 
\end{tabular} \caption{Relaxation energy per site for TFI models of several sizes with $h=0.5,1,1.5$
and clusters of fixed size $2\times1$.}
\label{tab:tfi_size_energies}
\end{table}

\appendix

\section{Fermions\label{sec:Fermions}}

In this appendix we present the relevant background on fermionic many-body
systems, following~\cite{LinLindsey2020}.

\subsection{Background}

Fermionic many-body problems in second quantization are specified
in terms of the creation operators $a_{1}^{\dagger},\ldots,a_{M}^{\dagger}$
and their Hermitian adjoints, the annihilation operators $a_{i}$,
which can be viewed as operators on a vector space of dimension $2^{M}$
called the Fock space $\mathcal{F}$. The key properties of these
operators are the canonical anticommutation relations 

\[
\{a_{i},a_{j}^{\dagger}\}=\delta_{ij},\quad\{a_{i},a_{j}\}=\{a_{i}^{\dagger},a_{j}^{\dagger}\}=0,
\]
 where $\{\,\cdot\,,\,\cdot\,\}$ denotes the anticommutator. In terms
of these operators we also define the number operators $\hat{n}_{i}:=a_{i}^{\dagger}a_{i}$,
and the total number operator by $\hat{N}:=\sum_{i=1}^{M}\hat{n}_{i}$.

One can identify the Fock space with a quantum spin-$\frac{1}{2}$
state space, i.e., identify $\mathcal{F}\simeq\bigotimes^{M}\mathbb{C}^{2}\simeq\mathbb{C}^{2^{M}}$,
via the correspondence 
\[
a_{i}^{\dagger}\rightsquigarrow\underbrace{\sigma^{z}\otimes\cdots\otimes\sigma^{z}}_{i-1\,\mathrm{factors}}\otimes\left(\begin{array}{cc}
0 & 0\\
1 & 0
\end{array}\right)\otimes I_{2}\otimes\cdots\otimes I_{2},
\]
 known as the Jordan-Wigner transformation (JWT). This transformation
is unnatural in the sense that it depends on the ordering of the states.
More precisely, permuting the states before the JWT is not equivalent
to permuting the tensor factors after the JWT. Moreover, a fermionic
operator such as $a_{i}^{\dagger}a_{j}$ involving only two sites
corresponds in general to a quantum spin operator involving potentially
many more sites; hence the pairwise fermionic Hamiltonians that we
consider below cannot be viewed in general as pairwise spin-$\frac{1}{2}$
Hamiltonians.

Next we define the notion of a pairwise fermionic Hamiltonian, and
then we provide some examples. Unfortunately we cannot simply treat
clusters of sites as `supersites' without breaking the fermionic structure,
so we approach the cluster framework directly, writing $\{1,\ldots,M\}$
as a disjoint union of clusters $\bigcup_{\gamma=1}^{N_{\mathrm{c}}}C_{\gamma}$
specified by the user. We let 
\[
\mathcal{A}:=\langle1,a_{1},\ldots,a_{M},a_{1}^{\dagger},\ldots,a_{M}^{\dagger}\rangle
\]
 denote the star-algebra generated by the creation and annihilation
operators subject to the canonical anticommutation relations $\{a_{i},a_{j}^{\dagger}\}=\delta_{ij}$,
$\{a_{i},a_{j}\}=\{a_{i}^{\dagger},a_{j}^{\dagger}\}=0$, and similarly,
we let 
\[
\mathcal{A}_{C}:=\left\langle \{1\}\cup\{a_{i},a_{i}^{\dagger}\,:\,i\in C\}\right\rangle 
\]
 denote the subalgebra corresponding to a subset $C\subset\{1,\ldots,M\}$.
Then we consider pairwise Hamiltonians $\hat{H}\in\mathcal{A}$ of
the form 
\[
\hat{H}=\sum_{\gamma}\hat{H}_{\gamma}+\sum_{\gamma<\delta}\hat{H}_{\gamma\delta},
\]
 where $\hat{H}_{\gamma}\in\mathcal{A}_{C_{\gamma}}$ and $\hat{H}_{\gamma\delta}\in\mathcal{A}_{C_{\gamma}\cup C_{\delta}}$
are Hermitian operators. We are interested in the ground-state energy
\[
E_{0}=\inf\left\{ \langle\psi\vert\hat{H}\vert\psi\rangle\,:\,\vert\psi\rangle\in\mathcal{F},\,\langle\psi\vert\psi\rangle=1\right\} .
\]
 For particle-number conserving Hamiltonians (i.e., Hamiltonians that
commute with $\hat{N}$), one may also consider the $N$-particle
ground-state energy defined as 
\[
E_{0}(N)=\inf\left\{ \langle\psi\vert\hat{H}\vert\psi\rangle\,:\,\vert\psi\rangle\in\mathcal{F},\,\langle\psi\vert\psi\rangle=1,\,\langle\psi\vert\hat{N}\vert\psi\rangle=N\right\} ,
\]
 though note that this is formally equivalent to the unconstrained
ground-state energy $E_{0}$ after subtracting $\mu\hat{N}$ from
the Hamiltonian, where the Lagrange multiplier $\mu$ is called the
chemical potential.

\subsection{Examples}

In this work we shall consider the half-filled lattice of spinless
fermions~\cite{SpinlessFerm} specified by the Hamiltonian 
\begin{equation}
\hat{H}=\sum_{i\sim j}\left[-a_{i}^{\dagger}a_{j}-a_{j}^{\dagger}a_{i}+U\left(\hat{n}_{i}-\frac{1}{2}\right)\left(\hat{n}_{j}-\frac{1}{2}\right)\right],\label{eq:spinlessFerm}
\end{equation}
 where $U$ is a scalar parameter (the `interaction strength') and
the notion of adjacency $i\sim j$ is defined relative to a graph
(usually a rectangular lattice) on the sites $\{1,\ldots,M\}$. This
operator is pairwise relative to any cluster decomposition. One can
also consider an analogous model with long-range Coulomb interaction
\begin{equation}
\hat{H}=\sum_{i\sim j}\left[-a_{i}^{\dagger}a_{j}-a_{j}^{\dagger}a_{i}\right]+U\sum_{i\neq j}\frac{1}{d(i,j)}\left(\hat{n}_{i}-\frac{1}{2}\right)\left(\hat{n}_{j}-\frac{1}{2}\right),\label{eq:spinlessFermLR}
\end{equation}
 where $d(i,j)$ is the Euclidean distance between sites $i$ and
$j$ on the lattice. 

As outlined in~\cite{LinLindsey2020}, one can also consider spinful systems such as the Hubbard model~\cite{Hubbard}, as well as quantum chemistry Hamiltonians arising from
electronic structure problems after a suitable choice of basis~\cite{SzaboOstlund1989}, e.g., the recently developed discontinuous Galerkin basis~\cite{DGBasis}.

\subsection{The two-marginal relaxation}

In order to realize the two-marginal relaxation as a concrete semidefinite
program, it is necessary to choose a JWT \emph{for each pair of clusters
}as follows. (Note that there will be no need to consider any global
JWT.) For $1\leq\gamma\leq N_{\mathrm{c}}$, let $L_{\gamma}:=\vert C_{\gamma}\vert$,
and let $\kappa_{\gamma}:C_{\gamma}\ra\{1,\ldots,L_{\gamma}\}$ be
a bijection specifying an ordering for the sites in the $\gamma$-th
cluster. For $1\leq\gamma<\delta\leq N_{\mathrm{c}}$, let $L_{\gamma\delta}:=\vert C_{\gamma}\vert+\vert C_{\delta}\vert$,
and let $\kappa_{\gamma\delta}:C_{\gamma}\cup C_{\delta}\ra\{1,\ldots,L_{\gamma\delta}\}$
be the bijection specifying an ordering for the sites in the $(\gamma,\delta$)-th
pair of clusters, uniquely specified by the conditions that $\kappa_{\gamma\delta}\vert_{C_{\gamma}}=\kappa_{\gamma}$,
$\kappa_{\gamma\delta}\vert_{C_{\delta}}=\kappa_{\delta}$, and $\kappa_{\gamma\delta}(C_{\gamma})<\kappa_{\gamma\delta}(C_{\delta})$.
These orderings fix algebra isomorphisms $\mathcal{J}_{\gamma}:\mathcal{A}_{C_{\gamma}}\ra\mathrm{End}\left(\bigotimes_{i=1}^{L_{\gamma}}\mathbb{C}^{2}\right)$
and $\mathcal{J}_{\gamma\delta}:\mathcal{A}_{C_{\gamma}\cup C_{\delta}}\ra\mathrm{End}\left(\bigotimes_{i=1}^{L_{\gamma\delta}}\mathbb{C}^{2}\right)$
via the appropriate JWTs. Then the two-marginal relaxation is given
concretely by 
\begin{eqnarray}
\underset{\{\rho_{\gamma}\},\,\{\rho_{\gamma\delta}\}_{\gamma<\delta}}{\mathrm{minimize}} &  & \sum_{\gamma}\Tr\left[\mathcal{J}_{\gamma}\left(\hat{H}_{\gamma}\right)\rho_{\gamma}\right]+\sum_{\gamma<\delta}\Tr\left[\mathcal{J}_{\gamma\delta}\left(\hat{H}_{\gamma\delta}\right)\rho_{\gamma\delta}\right],\label{eq:sdpFerm}\\
\mbox{subject to} &  & \rho_{\gamma\delta}\succeq0,\quad1\leq\gamma<\delta\leq N_{\mathrm{c}},\nonumber \\
 &  & \rho_{\gamma}=\Tr_{\kappa_{\gamma\delta}(C_{\delta})}[\rho_{\gamma\delta}],\quad\rho_{\delta}=\Tr_{\kappa_{\gamma\delta}(C_{\gamma})}[\rho_{\gamma\delta}],\quad1\leq\gamma<\delta\leq N_{\mathrm{c}},\nonumber \\
 &  & \Tr[\rho_{\gamma}]=1,\quad\gamma=1,\ldots,N_{\mathrm{c}},\nonumber \\
 &  & G\left[\{\rho_{\gamma}\},\{\rho_{\gamma\delta}\}_{\gamma<\delta}\right]\succeq0,\nonumber 
\end{eqnarray}
 where $G$ is specified blockwise subordinate to a collection $\{\hat{A}_{\gamma,\alpha}\,:\,\alpha=1,\ldots,n_{\gamma}\}\subset\mathcal{A}_{C_{\gamma}}$
of operators on each cluster via 
\[
\left(G_{\gamma\delta}\right)_{\alpha\beta}=\begin{cases}
\Tr\left(\left[\mathcal{J}_{\gamma}\left(\hat{A}_{\gamma,\alpha}\right)\right]^{\dagger}\left[\mathcal{J}_{\gamma}\left(\hat{A}_{\gamma,\beta}\right)\right]\rho_{\gamma}\right), & \gamma=\delta\\
\Tr\left(\left[\mathcal{J}_{\gamma\delta}\left(\hat{A}_{\gamma,\alpha}\right)\right]^{\dagger}\left[\mathcal{J}_{\gamma\delta}\left(\hat{A}_{\delta,\beta}\right)\right]\rho_{\gamma\delta}\right), & \gamma\neq\delta.
\end{cases}
\]
 Hence after the appropriate Jordan-Wignerized operators are formed,
(\ref{eq:sdpFerm}) is of identical form to (\ref{eq:sdpObj}).

\bibliography{varembed}
\bibliographystyle{siam}

\end{document}